\documentclass[10pt,11pt]{article}%
\usepackage[dvips]{graphicx}
\usepackage{amsmath}
\usepackage{color}
\usepackage{amsfonts}
\usepackage{amssymb}%
\providecommand{\U}[1]{\protect\rule{.1in}{.1in}}
\providecommand{\U}[1]{\protect\rule{.1in}{.1in}}
\providecommand{\U}[1]{\protect\rule{.1in}{.1in}}
\providecommand{\U}[1]{\protect\rule{.1in}{.1in}}
\providecommand{\U}[1]{\protect\rule{.1in}{.1in}}
\newtheorem{theorem}{Theorem}[subsection]

\newtheorem{condition}[theorem]{Condition}

\newtheorem{corollary}[theorem]{Corollary}

\newtheorem{example}[theorem]{Example}
\newtheorem{examples}[theorem]{Examples}

\newtheorem{lemma}[theorem]{Lemma}

\newtheorem{proposition}[theorem]{Proposition}

\newenvironment{proof}[1][Proof]{\textbf{#1.} }{\ \rule{0.5em}{0.5em}}

\def \Z{\mathbb{Z}}

\begin{document}

\title{Conditioned one-way simple random walk and representation theory}
\date{14/02/12}
\author{C\'{e}dric Lecouvey, Emmanuel Lesigne and Marc Peign\'{e}}
\maketitle

\begin{abstract}
We call one-way simple random walk a random walk in the quadrant
$\Z_+^n$ whose increments belong to the canonical base. In relation
with representation theory of Lie algebras and superalgebras, we
describe the law of such a random walk conditioned to stay in a closed
octant, a semi-open octant or other types of semi-groups. The combinatorial representation theory of these algebras allows us to describe a
generalized Pitman transformation which realizes the conditioning on
the set of paths of the walk. We pursue here in a direction initiated by O'Connell and his coauthors (\cite{OC1}, \cite{OC2}, \cite{BBOC1}), and also developed in (\cite{LLP}. Our work relies on crystal bases
theory and insertion schemes on tableaux described by Kashiwara and
his coauthors in \cite{BKK} and, very recently, in \cite {GHSKM}.
\end{abstract}


\section{Introduction}
Let $B=(\varepsilon_{1},\ldots,\varepsilon_{n})$ be the standard basis of
$\mathbb{R}^{n}$.\ The one-way simple walk is defined as the random walk
$\mathcal{W}=(\mathcal{W}_{\ell}=X_{1}+\cdots+X_{\ell})_{\ell\geq1}$ where $(X_{k}
)_{k\geq1}$ is a sequence of independent and identically distributed random
variables with values in the base $B$ and with common mean  vector $\mathbf{m}${.}  In this paper, we
generalize some results due to O'Connell \cite{OC1},\cite{OC2} given the law
of the random walk $\mathcal{W}$ conditioned
 never to exit the cone $\mathcal C^\emptyset=\{x=(x_{1},\ldots,x_{n})\in\mathbb{R}
^{n}\mid x_{1}\geq\cdots\geq x_{n}\geq0\}$\footnote{We will discuss in parallel the cases of Lie algebras $\mathfrak{gl}(n),\mathfrak{gl}(m,n)$ and $\mathfrak{q}(n)$, and an exponent $\emptyset$, $h$ or $s$ will refer respectively to each of these three choices. See $\S\,$\ref{Wsc}.}. This is achieved in \cite{OC1} by
considering first a natural transformation $\mathfrak{P}$ which associates to
any path with steps in $B$ a path in the cone $\mathcal C^\emptyset$, next by checking
that the image of the random walk $\mathcal{W}$ by this transformation
is a Markov chain and finally by establishing that this Markov chain has the
same law as $\mathcal{W}$ conditioned  never to  exit $\mathcal C^\emptyset$.\ Observe
that $\mathcal C^\emptyset\cap\mathbb{Z}^{n}$ is the set of partitions $\mathcal{P}
^{\emptyset}=\{\lambda=(\lambda_{1},\ldots,\lambda_{n})\in \mathbb Z^{n}\mid\lambda
_{1}\geq\cdots\geq\lambda_{n}\geq0)\}$. The entries of the transition matrix
of the Markov chain so obtained are indexed by the pairs of such partitions.
Moreover, they can be expressed as quotients of Schur functions (the Weyl
characters of $\mathfrak{sl}_{n}$) with variables specialized to the
coordinates of $\mathbf{m}$. The transformation $\mathfrak{P}$ is based on the
Robinson-Schensted correspondence which maps the words on the ordered alphabet
$\mathcal{A}_{n}=\{1<\cdots<n\}$ (regarded as finite paths in $\mathbb{Z}^{n}
$) on pairs of semistandard tableaux. Namely, the map $\mathfrak{P}$ associates to any
trajectory in $\mathbb{Z}^{n}$ of length $\ell$ (encoded by a word $w$ of
length $\ell$) the sequence of shapes of the tableaux obtained by applying
the Robinson-Schensted procedure to $w$; this sequence is then regarded as a
trajectory in $\mathcal C^\emptyset$ of length $\ell$. It is proved in
\cite{OC1} that, for $n=2$, the transformation $\mathfrak{P}$ coincides with
the usual Pitman transform on the line.

One can introduce a similar transformation $\mathfrak{P}$ for a wide class of
random walks {$(X_{1}+\cdots+X_{\ell})_{\ell\geq1}$ } for
which the variables $X_{k}$ take values in the set of weights of a fixed
representation of a simple Lie algebra $\mathfrak{g}(n)$ over $\mathbb{C}$.\ This
was done in \cite{BBOC1} in the case of equidistributed random variables
$X_{k}$ and in \cite{LLP} in general. The transformation $\mathfrak{P}$ is
then defined by using Kashiwara's crystal basis theory \cite{Kashi} (or
equivalently in terms of the Littelmann path model). We also obtain in
\cite{LLP} the law of the random walk   conditioned to never exit the cone $\mathcal C^\emptyset$ of dominant
integral weights for $\mathfrak{g}(n)$ under the crucial assumption (also
required in \cite{OC1} and \cite{OC2}) that $\mathbf{m}=E(X_{k})$ belongs to
 the interior  $C^\circ$ of $\mathcal C^\emptyset$. The transition matrix obtained has a
simple expression in terms of the Weyl characters of the irreducible
representations of $\mathfrak{g}(n)$. It is also worth mentioning that the
interaction between random walk problems and representation theory permits to
derive new results in both theories. In particular, the purely probabilistic
theorems established in \cite{LLP} also yield the asymptotic behavior of
certain tensor product multiplicities which seems very difficult to reach by
purely algebraic arguments.

It is then a natural question to try to extend the results of \cite{OC1} and
\cite{LLP} to other conditionings. In this paper, we generalize the results of
\cite{OC1} to the one-way simple walk conditioned to stay in discrete cones or
semigroups appearing naturally in the representation theories of the Lie
superalgebras $\mathfrak{gl}(m,n)$ and $\mathfrak{q}(n)$. 

For the superalgebra $\mathfrak{gl}(m,n)$, we give the law of the random walk $\mathcal{W}$
conditioned to stay in the semigroup $\mathcal{C}^h$ of $\mathbb{R}^{m+n}$ defined in $\S \, $\ref{tmP}; the subsemigroup $\mathcal P^{h}=\mathcal{C}^h\cap\Z^{m+n}$ is 
 parametrized by the hook partitions; the drift $\mathbf m$ of the random walk is supposed belonging to the interior of $\mathcal{C}^h$. 
 
 For the superalgebra
$\mathfrak{q}(n)$, we give the law of the random walk $\mathcal{W}$ with drift in $C^\circ$ conditioned to never exit
the subcone $\mathcal C^s$ of $\mathcal C^\emptyset$ of vectors whose nonzero
coordinates are distinct. 

In both cases, this is also
achieved by introducing a Pitman type transform $\mathfrak{P}$. There are
nevertheless important differences with the Pitman transforms used in
\cite{OC1},\cite{BBOC1} and \cite{LLP}; indeed, these transforms can be obtained by
interpreting each path as a vertex in a Kashiwara crystal and then by applying
raising crystal operators until the highest weight (source) vertex is
reached. A contrario, the situation is more complicated for $\mathfrak{gl}(m,n)$ and
$\mathfrak{q}(n)$ since the associated crystals may admit fewer highest weight
vertices so that a relevant notion of Pitman transform cannot be defined by
using only raising operators. Let us also observe that the Pitman transforms
corresponding to $\mathfrak{gl}(m,n)$ and $\mathfrak{q}(n)$ do not fix the
paths contained in $\mathcal{C}^h$ and $\mathcal{C}^s$ (resp.) whereas the path contained in
$\mathcal C^\emptyset$ are fixed by the Pitman transform corresponding to $\mathfrak{gl}(n)$.

To overcome the complications due to the existence of fewer highest
weight vertices in crystals, we will define the transformation $\mathfrak{P}$
by using analogues of the Robinson-Schensted insertion procedure on tableaux
introduced in \cite{BKK} and very recently in \cite{GHSKM}.\ In fact the
complete crystal basis theories for the superalgebras $\mathfrak{gl}(m,n)$ and
$\mathfrak{q}(n)$ developed in \cite{BKK} and \cite{GHSKM} are crucial
ingredients of the paper.

Entries of the transition matrix of the one-way simple walk $\mathcal{W}$
conditioned to stay in $C$ are indexed by the pairs $(\lambda,\mu
)\in\mathcal{P}^{s}$ of strict partitions. We prove that these entries can be
expressed in terms of $P$-Schur functions (i.e. the characters of some
natural irreducible representations of $\mathfrak{q}(n)$). In addition to the
combinatorial representation theory of the Lie algebra $\mathfrak{q}(n)$, our
proof requires a quotient local limit theorem for random walks conditioned to
stay in cones, established in \cite{LLP}.\ Our method thus differs from that
used in \cite{OC1} and \cite{OC2}.\ In particular, we obtain the asymptotic
behavior of the tensor multiplicities corresponding to the vector
representation as a consequence of our main results whereas it was used as a
key ingredient in O'Connell's work.

For the conditioning of the walk $\mathcal{W}$ to stay in
$\mathcal{C}$, we use the representation theory of the Lie superalgebra
$\mathfrak{gl}(m,n)$. Entries of the transition matrix are expressed in
terms of super Schur functions, that are characters of the irreducible
representations appearing in the tensor powers of the vector representation of
$\mathfrak{gl}(m,n)$.\ We derive in particular the asymptotic behavior of the
corresponding tensor multiplicities. Observe this requires to extend the
quotient local limit theorem of \cite{LLP} to the case of random walks
conditioned to stay in semigroups since $\mathcal{C}$ is not a cone.

\bigskip

We will study simultaneously the three conditionings of the one-way simple walk in
$\mathcal C^\emptyset,C$ and $\mathcal{C}$ by using representation theory of
$\mathfrak{gl}(n),\mathfrak{gl}(m,n)$ and $\mathfrak{q}(n)$,
respectively.\ In particular, we obtain a different proof of the results of
O'Connell avoiding delicate determinantal computations required to obtain the
asymptotic behavior of some tensor multiplicities in the case of
$\mathfrak{gl}(n)$. We introduce the generalized Pitman transformations by using
insertion procedures on tableaux and Robinson-Schensted correspondences (as in
\cite{OC1}) rather than the crystal basis theory (as in \cite{BBOC1} and
\cite{LLP}). We thus avoid technical difficulties inherent to the crystal
basis theory of the superalgebras $\mathfrak{gl}(m,n)$ and $\mathfrak{q}(n)$.
This also gives a different light on the Pitman transform(s) and a construction
step by step.

The paper is organized as follows. Sections \ref{Sec_Markov} and \ref{Sec_RT}
are respectively devoted to some background on Markov chains and
representation theory.\ In particular we establish Theorem \ref{RQT}, the
local limit theorem for a random walk conditioned to stay in a semigroup that
we need in the sequel.\ In Section \ref{sec_ComTab}, we recall the
combinatorial representation theory of the algebras $\mathfrak{gl}
(n),\mathfrak{gl}(m,n)$ and $\mathfrak{q}(n)$.\ In particular, we review the
relevant notions of Robinson-Schensted correspondences. Section
\ref{Sec_Pitman} \linebreak introduces the generalized Pitman transform. Without any
hypothesis on the drift, we
show that it maps the one-way simple walk onto a Markov chain whose transition matrix is computed. The
main result of the paper (Theorem \ref{Th_coincide}) giving the law of the
conditioned one-way simple walk with suitable drift is stated in Section
\ref{Sec_restric}.  Appendix is devoted to the proof of Proposition
\ref{Cor_RSK} for which complements on crystal basis theory for superalgebras
are required.

\section{Markov chains}

\label{Sec_Markov}We now recall the background on Markov chains and their
conditioning that we use in the sequel.\ 

\subsection{Markov chains and conditioning}

\label{mcc}

Consider a probability space $(\Omega,\mathcal{T},\mathbb{P})$ and a countable
set $M$. Let $Y=(Y_{\ell})_{\ell\geq1}$ be a sequence of random variables
defined on $\Omega$ with values in $M$. The sequence $Y$ is a Markov chain
when
\[
\mathbb{P}[Y_{\ell+1}=y_{\ell+1}\mid Y_{\ell}=y_{\ell},\ldots,Y_{1}
=y_{1}]=\mathbb{P}[Y_{\ell+1}=y_{\ell+1}\mid Y_{\ell}=y_{\ell}]
\]
for any $\ell\geq1$ and any $y_{1},\ldots,y_{\ell},y_{\ell+1}\in M$. The
Markov chains considered in the sequel will also be assumed time homogeneous,
that is $\mathbb{P}[Y_{\ell+1}=y_{\ell+1}\mid Y_{\ell}=y_{\ell}]=\mathbb{P}
[Y_{\ell}=y_{\ell+1}\mid Y_{\ell-1}=y_{\ell}]$ for any $\ell\geq2$.\ For all
$x,y$ in $M$, the transition probability from $x$ to $y$ is then defined by
\[
\Pi(x,y)=\mathbb{P}[Y_{\ell+1}=y\mid Y_{\ell}=x]
\]
and we refer to $\Pi$ as the transition matrix of the Markov chain $Y$. The
distribution of $Y_{1}$ is called the initial distribution of the chain $Y$.
It is well known that the initial distribution and the transition probability
determine the law of the Markov chain and that given a probability
distribution and a transition matrix on $M$, there exists an associated Markov
chain.\newline

Let $Y$ be a Markov chain on $(\Omega,\mathcal{T},\mathbb{P})$, whose initial
distribution has full support, i.e. $\mathbb{P}[Y_{1}=x]>0$ for any $x\in M$.
Let $\mathcal{C}$ be a nonempty subset of $M$ and consider the event
$S=[Y_{\ell}\in\mathcal{C}$ for any $\ell\geq1]$. Assume that $\mathbb{P}
[S\mid Y_{1}=\lambda]>0$ for all $\lambda\in\mathcal{C}$. This implies that
$\mathbb{P}[S]>0$, and we can consider the conditional probability
$\mathbb{Q}$ relative to this event: $\mathbb{Q}[\cdot]=\mathbb{P}[\cdot|S]$.

It is easy to verify that, under this new probability $\mathbb{Q}$, the
sequence $(Y_{\ell})$ is still a Markov chain, with values in $\mathcal{C}$,
and with transitions probabilities given by
\begin{equation}
\mathbb{Q}[Y_{\ell+1}=\lambda\mid Y_{\ell}=\mu]=\mathbb{P}[Y_{\ell+1}
=\lambda\mid Y_{\ell}=\mu]\ \frac{\mathbb{P}[S\mid Y_{1}=\lambda]}
{\mathbb{P}[S\mid Y_{1}=\mu]}. \label{reco}
\end{equation}
We will denote by $Y^{\mathcal{C}}$ this Markov chain and by $\Pi
^{\mathcal{C}}$ the restriction of the transition matrix $\Pi$ to the entries
which belong to $\mathcal{C}$ (in other words $\displaystyle
\Pi^{\mathcal{C}}=\left(  \Pi(\lambda,\mu)\right)  _{\lambda,\mu\in
\mathcal{C}}).$

\subsection{Doob $h$-transform}

\label{sub_sec_Doobh}

A \emph{substochastic matrix} on the countable set $M$ is a map $\Pi:M\times
M\rightarrow\lbrack0,1]$ such that $\sum_{y\in M}\Pi(x,y)\leq1$ for any $x\in
M.\;$If $\Pi,\Pi^{\prime}$ are substochastic matrices on $M$, we define their
product $\Pi\times\Pi^{\prime}$ as the substochastic matrix given by the
ordinary product of matrices:
\[
\Pi\times\Pi^{\prime}(x,y)=\sum_{z\in M}\Pi(x,z)\Pi^{\prime}(z,y).
\]

The matrix $\Pi^{\mathcal{C}}$ defined in the previous subsection is an
example of substochastic matrix.\newline

A function $h:M\rightarrow\mathbb{R}$ is \emph{harmonic} for the substochastic
transition matrix $\Pi$ when we have $\sum_{y\in M}\Pi(x,y)h(y)=h(x)$ for any
$x\in M$. Consider a strictly positive harmonic function $h$. We can then
define the Doob transform of $\Pi$ by $h$ (also called the $h$-transform of
$\Pi$) setting
\[
\Pi_{h}(x,y)=\frac{h(y)}{h(x)}\Pi(x,y).
\]
We then have $\sum_{y\in M}\Pi_{h}(x,y)=1$ for any $x\in M.\;$Thus $\Pi_{h}$
can be interpreted as the transition matrix for a certain Markov chain.

An example is given in the second part of the previous subsection (see formula
(\ref{reco})): the state space is now $\mathcal{C}$, the substochastic matrix
is $\Pi^{\mathcal{C}}$ and the harmonic function is $h_{\mathcal{C}}
(\lambda):=\mathbb{P}[S\mid Y_{1}=\lambda]$; the transition matrix
$\Pi_{h_{\mathcal{C}}}^{\mathcal{C}}$ is the transition matrix of the Markov
chain $Y^{C}$.

\subsection{Green function and Martin kernel}

Let $\Pi$ be a substochastic matrix on the set $M$. Its Green function is
defined as the series
\[
\Gamma(x,y)=\sum_{\ell\geq0}\Pi^{\ell}(x,y).
\]
If $\Pi$ is the transition matrix of a Markov chain, $\Gamma(x,y)$ is the
expected value of the number of passage at $y$ of the Markov chain starting at
$x$.

Assume that there exists $x^{\ast}$ in $M$ such that $0<\Gamma(x^{\ast
},y)<\infty$ for any $y\in M$. Fix such a point $x^{\ast}$.\ The Martin kernel
associated to $\Pi$ (with reference point $x^{\ast}$) is then defined by
\[
K(x,y)=\frac{\Gamma(x,y)}{\Gamma(x^{\ast},y)}.
\]
Consider a positive harmonic function $h$ and denote by $\Pi_{h}$ the
$h$-transform of $\Pi$. Consider the Markov chain $Y^{h}=\left(  Y_{\ell}
^{h}\right)  _{\ell\geq1}$ starting at $x^{\ast}$ and whose transition matrix
is $\Pi_{h}$. The following theorem is due to Doob. We gave in \cite{LLP} a
detailed proof which seems not very accessible in the literature.

\begin{theorem}
\label{Th_Doob}(Doob) Assume that there exists a function $f:M\rightarrow
\mathbb{R}$ such that for all $x\in M$, $\lim_{\ell\rightarrow+\infty
}K(x,Y_{\ell}^{h})=f(x)$ almost surely. Then there exists a positive real
constant $c$ such that $f=ch$.
\end{theorem}

\subsection{Quotient local limit theorem for a random walk in a semigroup}

\label{stayC} We now state some results on random walks similar to those
established in \cite{LLP}. However, the notion of random walk in a cone which
appears in \cite{LLP} will be replaced by the notion of random walk in a semigroup.

Let $\mathcal{C}$ be a subset of $\mathbb{R}^{n}$, stable under addition. We
suppose that its interior $\mathcal{C}^{\circ}$ is non empty. We denote by
$\mathcal{C}_{c}$ the cone of $\mathbb{R}^{n}$ generated by $\mathcal{C}$. We
will make the following assumption :

\begin{itemize}
\item[(h1)] $\mathcal{C}^{\circ}+\mathcal{C}_{c}\subset\mathcal{C}$.
\end{itemize}

(It would be sufficient for our purpose to suppose that for all $x\in
\mathcal{C}^{\circ}$ with $\|x\|$ large enough, and for all $y\in
\mathcal{C}_{c}$, we have $x+y\in\mathcal{C}$; but the examples that we will
consider satisfy (h1).)\newline

Here are two examples of additive subsemigroups of $\mathbb{R}^{n}$ that
will appear in the sequel of this article.

\begin{itemize}
\item[(e1)] In $\mathbb{R}^{n}$, any convex cone $\mathcal{C}$ with non empty
interior satisfies (h1).

\item[(e2)] Let $p,q$ be two positive integers, and $n=p+q$. If we define
\begin{multline*}
\mathcal{C}=\big\{(x_{\overline p},x_{\overline{p-1}},\ldots,x_{\overline1},x_{1},x_{2},\ldots,x_{q}
)\in\mathbb{R}^{p+q}\\
\mid x_{\overline p}\geq x_{\overline{p-1}}\geq\ldots\geq x_{\overline1}\geq0,\ x_{1}\geq x_{2}\geq
\ldots\geq x_{q}\geq0,\text{and }\forall i>x_{-1},x_{i}=0\big\},
\end{multline*}
then $\mathcal{C}$ is a semigroup. Moreover we have
\[
\mathcal{C}^{\circ}=\big\{x_{\overline p}>x_{\overline{p-1}}>\ldots>x_{\overline1}>q,\ x_{1}>x_{2}
>\ldots>x_{q}>0\big\},
\]
\[
\mathcal{C}_{c}\subset\big\{x_{\overline p}\geq x_{\overline{p-1}}\geq\ldots\geq x_{\overline1}
\geq0,\ x_{1}\geq x_{2}\geq\ldots\geq x_{q}\geq0\big\},
\]
and (h1) is satisfied.
\end{itemize}

\begin{lemma}
\label{sg-c}Denote by $\mathcal{C}_{c}^{\circ}$ the cone generated by
$\mathcal{C}^{\circ}$. For all compact set $K\subset\mathcal{C}_{c}^{\circ}$
and all large enough $t>0$, we have $tK\subset\mathcal{C}^{\circ}$.
\end{lemma}

A proof of Lemma \ref{sg-c} is given at the end of $\S\,$\ref{stayC}
.\newline

\medskip

Let $(X_{\ell})_{{\ell}\geq1}$ be a sequence of independent and identically
distributed random variables defined on a probability space $(\Omega
,{\mathcal{T}},\mathbb{P})$ and with values in the Euclidean space
${\mathbb{R}}^{n}$. We assume that these variables have a moment of order 1
and denote by $\mathbf{m}$ their common mean vector. Let us denote by $(S_{\ell
})_{{\ell}\geq0}$ the associated random walk defined by $S_{0}=0$ and
$S_{\ell}:=X_{1}+\cdots+X_{\ell}$ for $\ell\geq1$.

We consider a semigroup $\mathcal{C}$ in $\mathbb{R}^{n}$, with interior
$\mathcal{C}^{\circ}\neq\emptyset$. In order to see enough paths of the random
walk staying in the semigroup we need some additional assumptions. We assume
that :

\begin{enumerate}
\item[(h2)] $\exists t>0,\ t\mathbf{m}\in\mathcal{C}^{\circ}$.

\item[(h3)] $\exists\ell_{0}>0,\ \mathbb{P}\left[  S_{1}\in\mathcal{C}
,S_{2}\in\mathcal{C},\ldots,S_{\ell_{0}-1}\in\mathcal{C},S_{\ell_{0}}
\in\mathcal{C}^{\circ}\right]  >0.$
\end{enumerate}

\begin{lemma}One gets
\label{psc}
\[
\mathbb{P}\left[  \forall{\ell}\geq1,\ S_{\ell}\in\mathcal{C}\right]  >0.
\]

\end{lemma}

\begin{proof}
By hypothesis (h3), one may fix $a$ in $\mathcal{C}^{\circ}$ such that, for
all $\varepsilon>0$,
\[
\mathbb{P}\left[  S_{1}\in\mathcal{C},\ldots,S_{\ell_{0}-1}\in\mathcal{C}
,S_{\ell_{0}}\in B(a,\varepsilon)\right]  >0.
\]

By the strong law of large numbers, the sequence $\left(  \frac{1}{\ell
}S_{\ell}\right)  _{\ell\geq1}$ converges almost surely to $\mathbf{m}$.
Therefore, thanks to Hypothesis (h2) and Lemma \ref{sg-c}, almost surely, one
gets $S_{\ell}\in\mathcal{C}^{\circ}$ for any large enough ${\ell}$, that is
\[
\lim_{L\rightarrow+\infty}\mathbb{P}\left[  \forall{\ell}\geq L,\ S_{\ell}
\in\mathcal{C}^{\circ}\right]  =1.
\]

By Lemma \ref{sg-c}, we know that, for any $x\in\mathbb{R}^{n}$, for all large
enough $k\in\mathbb{N}$ one gets $x+ka\in\mathcal{C}^{\circ}$. Thus
\[
\left[  \forall{\ell}\geq L,\ S_{\ell}\in\mathcal{C}^{\circ}\right]
=\bigcup_{k\geq0}\left[  (\forall{\ell}\geq L,\ S_{\ell}\in\mathcal{C}^{\circ
})\text{ and }(\forall{\ell}<L,\ S_{\ell}+ka\in\mathcal{C}^{\circ})\right]
\]
and therefore $\displaystyle\left[  \forall{\ell}\geq L,\ S_{\ell}
\in\mathcal{C}^{\circ}\right]  \subset\bigcup_{k\geq0}\left[  \forall{\ell
},\ S_{\ell}+ka\in\mathcal{C}^{\circ}\right]  $ since $\mathcal{C}^{\circ}$ is
a semigroup. Hence, there exists $k\geq0$ such that $\mathbb{P}\left[
\forall{\ell},\ S_{\ell}+ka\in\mathcal{C}^{\circ}\right]  >0.$ Fix such a $k$
and $\delta>0$ such that $B(a,\delta)\subset\mathcal{C}^{\circ}$.

We consider now $k$ independent repetitions of the event
\[
A_{1}:=\left[  S_{1}\in\mathcal{C},\ldots,S_{\ell_{0}-1}\in\mathcal{C}
,\ S_{\ell_{0}}\in B(a,\frac{\delta}{k+1})\right]  ,
\]
that is we consider the events $A_{j}$, $1\leq j\leq k+1$, defined by
\[
A_{j}:=\left[  S_{(j-1)\ell_{0}+1}-S_{(j-1)\ell_{0}}\in\mathcal{C}
,\ldots,S_{(j-1)\ell_{0}+\ell_{0}-1}-S_{(j-1)\ell_{0}}\in\mathcal{C}
,\ S_{j\ell_{0}}-S_{(j-1)\ell_{0}}\in B(a,\frac{\delta}{k+1})\right]  .
\]
In the following claim, we use the semigroup property: if the events
$A_{1},\ldots,A_{k+1}$ are simultaneously realized, then
\[
S_{1}\in\mathcal{C},\ldots,S_{(k+1)\ell_{0}}\in\mathcal{C}\text{ and }\Vert
S_{(k+1)\ell_{0}}-ka-a\Vert\leq\sum_{j=1}^{k+1}\Vert S_{j\ell_{0}
}-S_{(j-1)\ell_{0}}-a\Vert\leq\delta,
\]
so $S_{(k+1)\ell_{0}}-ka\in\mathcal{C}^{\circ}$. We thus have

\begin{multline*}
\mathbb{P}\left[  \forall{\ell}\geq1,\ S_{\ell}\in\mathcal{C}\right]
\geq\mathbb{P}\left[  \left(  \cap_{j=1}^{k+1}A_{j}\right)  \text{ and
}\left(  \forall{\ell}>(k+1)\ell_{0},\ S_{\ell}-S_{(k+1)\ell_{0}}
+ka\in\mathcal{C}^{\circ}\right)  \right] \\
=\prod_{j=1}^{k+1}\mathbb{P}[A_{j}]\ \times\ \mathbb{P}\left[  \forall{\ell
}>0,\ S_{\ell}+ka\in\mathcal{C}^{\circ}\right]  =\left(  \mathbb{P}
[A_{1}]\right)  ^{k+1}\ \times\ \mathbb{P}\left[  \forall{\ell}>0,\ S_{\ell
}+ka\in\mathcal{C}^{\circ}\right]  >0.
\end{multline*}

\end{proof}

The quotient local limit theorem stated in \cite{LLP} can be extended to our
situation. We limit our study to random walks in the discrete lattice
$\mathbb{Z}^{n}$, and we have to make an aperiodicity hypothesis: we assume
that the support $S_{\mu}$ of the law $\mu$ of the random variables $X_{\ell}$
is a subset of $\mathbb{Z}^{n}$ and that $S_{\mu}$ is not contained in a
coset of a proper subgroup of $\mathbb{Z}^{n}$.

\begin{theorem}
\label{RQT} Assume that the random variables $X_{\ell}$ are almost surely
bounded. Let $\mathcal{C}$ be an additive subsemigroup of $\mathbb{R}^{n}$
satisfying hypothesis (h1), (h2) and (h3). Let $(g_{\ell}),(h_{\ell})$ be two
sequences in $\mathbb{Z}^{n}$ and $\alpha<2/3$ such that $\lim{\ell}^{-\alpha
}\Vert g_{\ell}-{\ell}\mathbf{m}\Vert=0$ and $\lim{\ell}^{-1/2}\Vert h_{\ell
}\Vert=0$. Then, when ${\ell}$ tends to infinity, we have
\[
\mathbb{P}\left[  S_{1}\in\mathcal{C},\ldots,S_{\ell}\in\mathcal{C},S_{\ell
}=g_{\ell}+h_{\ell}\right]  \sim\mathbb{P}\left[  S_{1}\in\mathcal{C}
,\ldots,S_{\ell}\in\mathcal{C},S_{\ell}=g_{\ell}\right]  .
\]

\end{theorem}

Some comments are necessary in order to justify this statement. Indeed, this
theorem is proved in \cite{LLP} in the case when $\mathcal{C}$ is a cone, but we relax this condition here. First of all, we claim that Lemma 4.4 from
\cite{LLP} still holds in our context of semigroups. Let us recall this statement.

\begin{lemma}
\label{lem-marc} Assume the random variables $X_{\ell}$ are almost surely
bounded. Let $\alpha\in]1/2,2/3[$. If the sequence $({\ell}^{-\alpha}\Vert
g_{\ell}-{\ell}\mathbf{m}\Vert)_{\ell}$ is bounded, then there exists $c>0$
such that, for all large enough $\ell$,
\begin{equation}
\mathbb{P}\left[ S_{1}\in\mathcal{C},\ldots,S_{\ell}\in\mathcal{C},S_{\ell
}=g_{\ell}\right]  \geq\exp\left(  -c\ell^{\alpha}\right)  . \label{garbit}
\end{equation}

\end{lemma}

\begin{proof}
Thanks to hypothesis (h2) we can fix $\delta>0$ such that the closure
$\overline{B}(\mathbf{m},\delta)$ of the open ball $B(\mathbf{m},\delta)$ is
contained in the cone $\mathcal{C}_{c}^{\circ}$ generated by $\mathcal{C}
^{\circ}$. Thanks to Lemma \ref{sg-c}, we know that $B(t\mathbf{m},t\delta
)\subset\mathcal{C}^{\circ}$ for all large enough $t$.

On the other hand, the law of large numbers tells us that $\displaystyle\lim
_{\ell\rightarrow+\infty}\mathbb{P}\left[ \left\Vert \frac{1}{\ell}S_{\ell
}-\mathbf{m}\right\Vert <\delta\right]  =1$. Combined with Lemma \ref{psc},
this implies that, for $\ell_{1}$ large enough,
\[
\mathbb{P}\left[  S_{1}\in\mathcal{C},\ldots,S_{\ell_{1}}\in\mathcal{C}
,S_{\ell_{1}}\in B(\ell_{1}\mathbf{m},\ell_{1}\delta)\right] >0.
\]

From the two preceding claims, we deduce the existence of an integer $\ell
_{1}>0$ and a point $y_{1}\in B(\ell_{1}\mathbf{m},\ell_{1}\delta
)\cap\mathbb{Z}^{n}$ such that
\[
\mathbb{P}\left[  S_{1}\in\mathcal{C},\ldots,S_{\ell_{1}}\in\mathcal{C}
,S_{\ell_{1}}=y_{1}\right]  >0\quad\mathrm{and}\quad B(t\mathbf{m},t\delta
)\subset\mathcal{C}^{\circ}\quad\mathrm{forall}\quad t\geq\ell_{1}.
\]
We fix such a pair $(\ell_{1},y_{1})$. Since $y_{1}\in\mathcal{C}^{\circ}$ we
have $y_{1}+\mathcal{C}_{c}\subset\mathcal{C}$, by hypothesis (h1). For
$\ell\geq\ell_{1}$, we have
\begin{multline*}
\mathbb{P}\left[  S_{1}\in\mathcal{C},\ldots,S_{\ell}\in\mathcal{C},S_{\ell
}=g_{\ell}\right] \\
\geq\mathbb{P}\big[S_{1}\in\mathcal{C},\ldots,S_{\ell_{1}}\in\mathcal{C}
,S_{\ell_{1}}=y_{1},S_{\ell_{1}+1}\in y_{1}+\mathcal{C}_{c},S_{\ell_{1}+2}\in
y_{1}+\mathcal{C}_{c},\ldots\\
\ldots,S_{\ell}\in y_{1}+\mathcal{C}_{c},S_{\ell}-S_{\ell_{1}}=g_{\ell}
-y_{1}\big]\\
=\mathbb{P}\left[  S_{1}\in\mathcal{C},\ldots,S_{\ell_{1}}\in\mathcal{C}
,S_{\ell_{1}}=y_{1}\right]  \times\mathbb{P}\left[ S_{1}\in\mathcal{C}
_{c},S_{2}\in\mathcal{C}_{c},\ldots,S_{\ell-\ell_{1}}\in\mathcal{C}
_{c},S_{\ell-\ell_{1}}=g_{\ell}-y_{1}\right]
\end{multline*}

Now, since the sequence $\left(  \ell^{-\alpha}\Vert g_{\ell+\ell_{1}}
-y_{1}-\ell\mathbf{m}\Vert\right)  $ is bounded, one may apply Lemma 4.4 of
\cite{LLP}, with the cone $\mathcal{C}_{c}$, which satisfies the required
properties with respect to the random walk; one gets
\[
\liminf_{\ell\rightarrow\infty}\left(  \mathbb{P}\left[ S_{1}\in
\mathcal{C}_{c},\ldots,S_{\ell-\ell_{1}}\in\mathcal{C}_{c},S_{\ell-\ell_{1}
}=g_{\ell}-y_{1}\right]  \right)  ^{{\ell}^{-\alpha}}>0,
\]
so that
\[
\liminf_{\ell\rightarrow\infty}\left(  \mathbb{P}\left[ S_{1}\in
\mathcal{C},\ldots,S_{\ell}\in\mathcal{C},S_{\ell}=g_{\ell}\right]\right)
^{{\ell}^{-\alpha}}>0.
\]

\end{proof}

\begin{proof}
[Proof of Theorem \ref{RQT}]We follow the lines of the proof of Theorem 4.3 in
\cite{LLP}. The claim \textquotedblleft For all $\ell\geq1$, $B(\ell
\mathbf{m},\ell\delta)\subset\mathcal{C}$\textquotedblright\ has to be
replaced by \textquotedblleft For all large enough $\ell$, $B(\ell
\mathbf{m},\ell\delta)\subset\mathcal{C}$\textquotedblright\ but this does not
disrupt the proof. In the proof of Theorem 4.3 in \cite{LLP}, we use Lemma
\ref{lem-marc} and the following claim: there exists $\varepsilon>0$ such that
when the walk goes out of $\mathcal{C}$, its distance to the point $k\mathbf m$ is
greater than $k\varepsilon$. As it is stated, this claim is not necessary true
in our present situation; fortunately this is certainly true for all large
enough $k$ and it is applied only to values of $k$ greater than $\ell^{\alpha
}$. For all the rest, the proof in \cite{LLP} can be followed line by line.
\end{proof}

\begin{proof}
[Proof of Lemma \ref{sg-c}]We have $K\subset\mathcal{C}_{c}^{\circ}$. For all
$x\in K$, there exists $t>0$ such that $tx\in\mathcal{C}^{\circ}$, and there
is a closed ball centered at $tx$ and contained in $\mathcal{C}^{\circ}$.
Coming back to $x$ we obtain the following : for all $x\in K$, there exist
$r>0$ and $t>0$ such that $t\overline{B}(x,r)\subset\mathcal{C}^{\circ}$. We
can cover $K$ by finitely many such balls $B(x,r)$, and it is sufficient to
prove the expected property for each of these balls.

We start from the fact that, for a particular $t=t_{0}>0$ we have
$t\overline{B}(x,r)\subset\mathcal{C}^{\circ}$, and we want to prove that it
is true for all large enough $t$.

We fix $\varepsilon>0$ such that $t_{0}\overline{B}(x,r+\varepsilon
)\subset\mathcal{C}^{\circ}$. Let $\vec{u}$ denote a unitary vector in
$\mathbb{R}^{n}$. For each $\vec{u}$ such that the half-line $\mathbb{R}
^{+}\vec{u}$ meets the ball $\overline{B}(x,r)$, we have
\begin{equation}
\mathbb{R}^{+}\vec{u}\cap t_{0}\overline{B}(x,r+\varepsilon)=[a(\vec
{u}),b(\vec{u})]\vec{u} \label{an1}
\end{equation}
where $a$ and $b$ are two continuous functions with $0<a<b$. By compacity, the
quotient $b/a$ stays greater than a number $\rho>1$ when the vector $\vec{u}$
varies as above. By the semigroup property of $\mathcal{C}^{\circ}$, for all
positive integer $n$,
\begin{equation}
\lbrack na(\vec{u}),nb(\vec{u})]\vec{u}\subset\mathcal{C}^{\circ}. \label{an2}
\end{equation}
Fix an integer $n_{0}>{\frac{1}{\rho-1}}.$ For all $n\geq n_{0}$, we have
$nb\geq(n+1)a$, hence
\begin{equation}
\cup_{n\geq n_{0}}[na,nb]=[n_{0}a,+\infty). \label{an3}
\end{equation}
From (\ref{an2}) and (\ref{an3}), we deduce that for all real $s\geq n_{0}$,
\[
s[a(\vec{u}),b(\vec{u})]\vec{u}\subset\mathcal{C}^{\circ}.
\]
With the help of (\ref{an1}), we conclude that for all real $s\geq n_{0}$,
\[
st_{0}\overline{B}(x,r)\subset\mathcal{C}^{\circ}.
\]
This concludes the proof.
\end{proof}


\section{Basics on representation theory}

\label{Sec_RT} We recall in the following paragraphs some classical material
on representation theory of classical Lie algebras and superalgebras needed in
the sequel.\ For a complete review, the reader is referred to \cite{Bour},
\cite{Hal} and \cite{Kac}.

\subsection{Weights and roots}

To the Lie algebra $\mathfrak{gl}(n)$ over $\mathbb{C}$ is associated its root
system. This root system is realized in an Euclidean space $\mathbb{R}^{n}$
with standard basis $B=(\varepsilon_{1},\ldots,\varepsilon_{n}).\;$The root
lattice of $\mathfrak{gl}(n)$ is the integral lattice $Q=\bigoplus_{i=1}
^{n-1}\mathbb{Z\alpha}_{i}$ where $\alpha_{i}=\varepsilon_{i}-\varepsilon
_{i+1}$ for $i=1,\ldots,n-1$.\ The weight lattice associated to $\mathfrak{gl}
(n)$ is the integral lattice $P=\mathbb{Z}^{n}.$ The cone of dominant positive
weights for $\mathfrak{gl}(n)$ is
\[
P_{+}=\{x=(x_{1},\ldots,x_{n})\in\mathbb{Z}^{n}\mid x_{1}\geq\cdots\geq
x_{n}\}\text{.}
\]
We also recall that the Weyl group of $\mathfrak{gl}(n)$ can be identified
with the symmetric group on $S_{n}.$ We write $\varepsilon$ for the signature
on the elements of $S_{n}$. The symmetric group $S_{n}$ acts on $\mathbb{Z}
^{n}$ by permutation on the coordinates. We write $\sigma_{0}$ the permutation
such that $\sigma_{0}(\beta_{1},\beta_{2},\ldots,\beta_{n})=(\beta_{n},\beta_{n-1}
,\ldots,\beta_{1})$ for any $(\beta_{1},\ldots,\beta_{n}
)\in\mathbb{Z}^{n}$. The Cartan Lie subalgebra $\mathfrak{h}$ of
$\mathfrak{gl}(n)$ is the subalgebra of diagonal matrices.\ The triangular
decomposition of $\mathfrak{gl}(n)=\mathfrak{gl}(n)_{+}\oplus\mathfrak{h}
\oplus\mathfrak{gl}(n)_{-}$ is the usual one obtained by considering strictly upper, diagonal and strictly lower matrices.

Now, we consider the Lie superalgebra $\mathfrak{gl}(m,n)$.\ It can be
regarded as the graded $\mathbb{Z}/2\mathbb{Z}$ algebra of the matrices of the
form
\[
\left(
\begin{array}
[c]{ll}
A & C\\
D & B
\end{array}
\right)  \quad A\in\mathfrak{gl}(m),B\in\mathfrak{gl}(n),C\in M(m,n),D\in
M(n,m)
\]
where $M(m,n)$ is the set of complex $m\times n$ rectangular matrices.\ It
decomposes as the sum of its even and odd parts $\mathfrak{gl}(m,n)=\mathfrak{gl}(m,n)
_{0}\oplus\mathfrak{gl}(m,n)_{1}$ where
\begin{align*}
\mathfrak{gl}(m,n)_{0}  &  =\left(
\begin{array}
[c]{ll}
A & 0\\
0 & B
\end{array}
\right)  \simeq\mathfrak{gl}(m)\oplus\mathfrak{gl}(n)\\
\text{ and} \qquad\mathfrak{gl}(m,n)_{1}  &  =\left(
\begin{array}
[c]{ll}
0 & C\\
D & 0
\end{array}
\right)  \simeq M(m,n)\oplus M(n,m).
\end{align*}
The ordinary Lie bracket is replaced by its super version, that is
$[X,Y]=XY-(-1)^{ij}YX$ for $X\in\mathfrak{gl}(n)_{i}$ and $Y\in\mathfrak{gl}(n)_{j}
$.\ Here $i$ and $j$ are regarded as elements of $\mathbb{Z}/2\mathbb{Z}
$.\ The Cartan subalgebra $\mathfrak{h},$ the Weyl group $W$ and the weight
lattice $P$ of $\mathfrak{gl}(m,n)$ coincide with those of the even part
$\mathfrak{gl}(n)_{0}$ and will be identified with $\mathbb{Z}^{n+m}$; in
particular $W=S_{m}\times S_{n}.$ In the sequel, it will be convenient to
write each weight $\beta\in\mathbb{Z}^{m+n}$ under the form $\beta
=(\beta_{\overline{m}},\ldots,\beta_{\overline{1}}\mid\beta_{1},\ldots
,\beta_{n})$. The set $P_{+}$ of positive dominant weights of $\mathfrak{gl}
(m,n)$ is also the same as the set of dominant weights of $\mathfrak{gl}(n)_{0}.$
We thus have $P_{+}=\{(\beta_{\overline{m}},\ldots,\beta_{\overline{1}}
\mid\beta_{1},\ldots,\beta_{n})\in\mathbb{Z}^{n+m}$ with $\beta_{\overline{m}
}\geq\cdots\geq\beta_{\overline{1}}$ and $\beta_{1}\geq\cdots\geq\beta_{n}\}$.
The superalgebra $\mathfrak{gl}(m,n)$ admits the triangular decomposition
$\mathfrak{gl}(m,n)=\mathfrak{gl}(m,n)_{+}\oplus\mathfrak{h}\oplus
\mathfrak{gl}(m,n)_{-}$ where $\mathfrak{gl}(m,n)_{+}$ and $\mathfrak{gl}
(m,n)_{-}$ are respectively the strictly upper and lower matrices in
$\mathfrak{gl}(m,n)$ and $\mathfrak{h}$ the subalgebra of diagonal matrices.

Here comes our third example. We denote by $\mathfrak{q}(n)$ the Lie
superalgebra of all matrices of the form
\[
\left(
\begin{array}
[c]{ll}
A & A^{\prime}\\
A^{\prime} & A
\end{array}
\right) \ , \quad A,A^{\prime}\in\mathfrak{gl}(n)\ ,
\]
endowed with the previous super Lie bracket.\ This superalgebra decomposes as
the sum $\mathfrak{q}(n)=\mathfrak{q}_{0}(n)\oplus\mathfrak{q}_{1}(n)$ of even
and odd parts with
\[
\mathfrak{q}_{0}(n)=\left(
\begin{array}
[c]{ll}
A & 0\\
0 & A
\end{array}
\right)  \simeq\mathfrak{gl}(n)\ \text{ and }\ \mathfrak{q}_{1}(n)=\left(
\begin{array}
[c]{ll}
0 & A^{\prime}\\
A^{\prime} & 0
\end{array}
\right)  .
\]
We denote by $e_{r,s},1\leq r,s\leq n$ (resp. $e_{r,s}^{\prime}$) the matrix
of $\mathfrak{q}_{0}(n)$ (resp. of $\mathfrak{q}_{1}(n)$) in which $A$ (resp.
$A^{\prime}$) has $(r,s)$-entry equal to $1$ and the other entries equal to
$0$. Let $\mathfrak{h}_{q}=\oplus_{r=1}^{n}\mathbb{C}e_{r,r}\oplus_{i=1}
^{n}\mathbb{C}e_{r,r}^{\prime}$ be the Cartan subalgebra of $\mathfrak{q}(n)$.
The superalgebra $\mathfrak{q}(n)$ admits the triangular decomposition
$\mathfrak{q}(n)=\mathfrak{q}_{+}(n)\oplus\mathfrak{h}_{q}\oplus
\mathfrak{q}_{-}(n)$ where $\mathfrak{q}_{+}(n)$ and $\mathfrak{q}_{-}(n)$ are
the subalgebras generated over $\mathbb{C}$ by $\{e_{r,s},e_{r,s}^{\prime}
\mid1\leq r<s\leq n\}$ and $\{e_{r,s},e_{r,s}^{\prime}\mid1\leq s<r\leq n\}.$
We also set $\mathfrak{h}=\oplus_{r=1}^{n}\mathbb{C}e_{r,r}.\;$The weight
lattice $P$ of $\mathfrak{q}(n)$ can be identified with $\mathbb{Z}^{n}$ and
its Weyl group with $S_{n}$. Both coincide with those of the even part
$\mathfrak{q}_{0}(n)$ and we set $P_{+}=\mathbb{Z}_{\geq0}^{n}.$

\subsection{Weight spaces and characters}\label{Wsc}

Assume $\mathfrak{g}=\mathfrak{gl}(n),\mathfrak{gl}(m,n)$ or $\mathfrak{q}
(n)$. For short, we set $N=n$ when $\mathfrak{g}=\mathfrak{gl}(n)$ or
$\mathfrak{q}(n)$ and $N=m+n$ when $\mathfrak{g}=\mathfrak{gl}(m,n)$. It will also be convenient to associate a symbol $\diamondsuit
=\emptyset,h,s$ to the objects attached to the algebras $\mathfrak{gl}
(n),\mathfrak{gl}(m,n),\mathfrak{q}(n),$ respectively.

In the sequel, we will only consider finite dimensional weight $\mathfrak{g}
$-modules.\ Such a module $M$ admits a decomposition in weight spaces
\[
M=\bigoplus_{\mu\in P}M_{\mu}\;,
\]
\[
M_{\mu}:=\{v\in M\mid h\cdot v=\mu(h)v\text{ for any }h\in\mathfrak{h}\}\;,
\]
where $P$ is embedded in the dual of $\mathfrak{h.}$ The space $M_{\mu}$ is
thus a $\mathfrak{h}$-module. If $M,M^{\prime}$ are finite-dimensional weight
$\mathfrak{g}$-modules and $\mu\in P$, we get $(M\oplus M^{\prime})_{\mu
}=M_{\mu}\oplus M_{\mu}^{\prime}$.\ In particular, the weight spaces
associated to any $\mathfrak{gl}(m,n)$ (resp.$\mathfrak{q}(n)$) module are
defined as the weight spaces of its restriction to $\mathfrak{gl}(m,n)_{0}$
(resp. $\mathfrak{q}(n)_{0}$). The character of $M$ is the Laurent polynomial
in $N$-variables
\[
\mathrm{char}(M)(x):=\sum_{\mu\in P}\dim(M_{\mu})x^{\mu}
\]
where $\dim(M_{\mu})$ is the dimension of the weight space $M_{\mu}$ and
$x^{\mu}$ is a formal exponential such that
\[
x^{\mu}=\left\{
\begin{array}
[c]{l}
x_{1}^{\mu_{1}}\cdots x_{n}^{\mu_{n}}\text{ for }\diamondsuit=\emptyset,s,\\
x_{\overline{m}}^{\mu_{\overline{m}}}\cdots x_{\overline{1}}^{\mu
_{\overline{1}}}x_{1}^{\mu_{1}}\cdots x_{n}^{\mu_{n}}\text{ for }
\diamondsuit=h.
\end{array}
\right.
\]
That is $x^{\mu}$ is determined by the coordinates of $\mu$ on the standard
basis $\{\varepsilon_{1},\ldots,\varepsilon_{N}\}.$ Here we identify as usual
the real form $\mathfrak{h}_{\mathbb{R}}^{\ast}$ of $\mathfrak{h}^{\ast}$ with
$\mathbb{R}^{N}$. For any $\sigma$ in the Weyl group, we have $\dim(M_{\mu
})=\dim(M_{\sigma(\mu)})$ so that $\mathrm{char}(M)(x)$ is a symmetric
polynomial for $\diamondsuit=\emptyset,s$ and is invariant under the action of
$S_{m}\times S_{n}$ for $\diamondsuit=h$.

A weight $\mathfrak{g}$-module $M$ is called a highest weight module with
highest weight $\lambda$ if $M$ is generated by $M_{\lambda}$ and
$\mathfrak{g}_{+}\cdot v=0$ for any $v$ in $M_{\lambda}$. To each dominant
weight $\lambda\in P_{+}$ corresponds a unique (up to isomorphism) irreducible
finite dimensional representation of $\mathfrak{g}$ of highest weight
$\lambda$. We denote it by $V^{\diamondsuit}(\lambda)$.\ We then write
\begin{equation}
\mathrm{char}(V^{\diamondsuit}(\lambda))(x):=\sum_{\mu\in P}K_{\lambda,\mu
}^{\diamondsuit}x^{\mu}. \label{Kostka}
\end{equation}
Thus $K_{\lambda,\mu}^{\diamondsuit}=\mathrm{dim(}V^{\diamondsuit}
(\lambda)_{\mu})$.

\subsection{Partitions and Young diagrams}

\label{subsec_YD}For any positive integer $k$, we denote by $\mathcal{P}_{k}
$ the set of partitions of length $k$. Recall that a partition $\lambda
\in\mathcal{P}_{k}$ is a $k$-tuple $\lambda=(\lambda_{1},\ldots,\lambda
_{k})\in\mathbb{Z}_{\geq0}^k$ such that $\lambda_{1}\geq\cdots\geq\lambda
_{k}$. We then set $\left\vert \lambda\right\vert =\lambda_{1}+\cdots
+\lambda_{k}$.\ The Young diagram $Y(\lambda)$ associated to $\lambda$ is the
juxtaposition of rows of lengths $\lambda_{1},\ldots,\lambda_{k}$ respectively, pictured
from top to bottom. Each row $i=1,\ldots,k$ is divided into $\lambda_{i}$
boxes and the rows are left justified. (See example below). The partition
$\lambda^{\prime}$ obtained by counting the number of boxes in each column of
$Y(\lambda)$ is the conjugate partition of $\lambda$.

$\bullet$ We write $\mathcal{P}^{\emptyset}=\mathcal{P}_{n}.\;$To each
partition $\lambda\in\mathcal{P}^{\emptyset}$, we associate the weight
$\pi(\lambda)=\sum_{i=1}^{n}\lambda_{i}\varepsilon_{i}$ which is dominant for
$\mathfrak{g}=\mathfrak{gl}(n)$.

One says that \textit{$Y(\lambda)$ is a $\emptyset$-diagram for $\mathfrak{gl}
(n)$.}

$\bullet$ Let $m,n$ be positive integers. We define the set $\mathcal{P}^{h}$
of \textit{hook partitions} as the set of partitions $\lambda$ of arbitrary
length such that $\lambda_{i}\leq n$ for any $i>m$.\ For any $\lambda
\in\mathcal{P}^{h}$, we denote by $\lambda^{(1)}\in\mathcal{P}_{m}$ the
partition corresponding to the Young diagram obtained by considering the $m$
longest rows of $\lambda$. We also denote by $\nu_{n}(\lambda)$ the partition
attached to the Young diagram obtained by deleting in the Young diagram of
$\lambda$ the boxes corresponding to $\lambda^{(1)}.$ By definition of the
hook partition $\lambda$, the conjugate partition of $\nu_{n}(\lambda)$
belongs to $\mathcal{P}_{n}$, it is denoted $\lambda^{(2)}$ and one gets
$\lambda^{(2)}=\nu_{n}(\lambda)^{\prime}$.\ We can write $\lambda
^{(1)}=(\lambda_{\overline{m}}^{(1)},\ldots,\lambda_{\overline{1}}^{(1)})$ and
$\lambda^{(2)}=(\lambda_{1}^{(2)},\ldots,\lambda_{n}^{(2)})$.\ We will then
write $\lambda=(\lambda^{(1)}\mid\lambda^{(2)})$ for short.\ This permits to
associate to $\lambda$ the weight $\pi(\lambda)=\sum_{i=1}^{m}\lambda
_{\overline{i}}^{(1)}\varepsilon_{i}+\sum_{j=1}^{n}\lambda_{i}^{(2)}
\varepsilon_{j+m}$, which is dominant for $\mathfrak{gl}(m,n)$.

One says that \textit{$Y(\lambda)$ is a $h$-diagram for $\mathfrak{gl}(m,n)$.}

$\bullet$ We similarly define the set $\mathcal{P}^{s}$ of \textit{strict
partitions} as the set of partitions $\lambda\in\mathcal{P}^{\emptyset}$ such
that $\lambda_{i+1}>0\Longrightarrow\lambda_{i}>\lambda_{i+1}$ for any
$i=1,\ldots,n-1.\;$We then define the shifted Young diagram $Y(\lambda)$
associated to $\lambda\in\mathcal{P}^{s}$ as the juxtaposition of rows of
lengths $\lambda_{1},\ldots,\lambda_{n}$ pictured from top to bottom. Each row
$i=1,\ldots,n$ is divided into $\lambda_{i}$ boxes but the $i$-th row is
shifted $i-1$ units to the right with respect to the top row. We also denote
by $\pi(\lambda)=\sum_{i=1}^{n}\lambda_{i}\varepsilon_{i}$ the dominant weight
of $\mathfrak{q}(n)$ associated to $\lambda$.

One says that \textit{$Y(\lambda)$ is a $s$-diagram for $\mathfrak{q}(n)$.}

\begin{example}
(1) The diagram $\displaystyle\quad
\begin{tabular}
[c]{|l|ll}\hline
&  & \multicolumn{1}{|l|}{}\\\hline
&  & \multicolumn{1}{|l|}{}\\\hline
&  & \multicolumn{1}{|l}{}\\\cline{1-2}
&  & \\\cline{1-1}
\end{tabular}
\quad$ is a $\emptyset$-diagram for $\mathfrak{gl}(4)$ with $\lambda
=(3,3,2,1).$

(2) The diagram $\displaystyle\quad
\begin{tabular}
[c]{|l|ll}\hline
&  & \multicolumn{1}{|l|}{}\\\hline
&  & \multicolumn{1}{|l|}{}\\\hline
&  & \multicolumn{1}{|l}{}\\\cline{1-2}
&  & \multicolumn{1}{|l}{}\\\cline{1-2}
&  & \multicolumn{1}{|l}{}\\\cline{1-2}
&  & \\\cline{1-1}
\end{tabular}
\quad$ is a $h$-diagram for $\mathfrak{gl}(2,2)$ with $\lambda
=(3,3,2,2,2,1).$

We have $\lambda^{(1)}=(3,3)$ and $\lambda^{(2)}=(4,3)$; as a element of
$\mathbb{Z}^{4}$, we thus have $\lambda=(3,3\mid4,3)$ and $\pi(\lambda
)=3\varepsilon_{\overline{2}}+3\varepsilon_{\overline{1}}+4\varepsilon
_{1}+3\varepsilon_{2}$

(3) The diagram $\displaystyle\quad
\begin{tabular}
[c]{ll|l|l}\hline
\multicolumn{1}{|l}{} & \multicolumn{1}{|l|}{} &  & \multicolumn{1}{|l|}{}
\\\hline
& \multicolumn{1}{|l|}{} &  & \\\cline{2-3}\cline{3-3}
&  &  & \\\cline{3-3}
\end{tabular}
\quad$ is a $s$-diagram for $\mathfrak{q}(3)$ with $\lambda=(4,2,1)$.
\end{example}

\noindent\textbf{Notation:} To simplify the notation, we shall identify in the
sequel the partition $\lambda$ with its associated dominant weight and simply
write $V^{\diamondsuit}(\lambda)$ for the highest weight module with highest
weight $\lambda$ rather than $V^{\diamondsuit}(\pi(\lambda)).$

\bigskip

\noindent\textbf{Remark: }For any $\diamondsuit=\{\emptyset,h,s\},$ the set
$\mathcal{P}^{\diamondsuit}$ is naturally associated to a subsemigroup $\mathcal C^\diamondsuit$ of the Euclidean space such that
$\mathcal{P}^{\diamondsuit}$ is the intersection of $\mathcal C^\diamondsuit$ with the integral lattice,
and $\mathcal C^\diamondsuit$ satisfies hypothesis (h1) of \S\,\ref{stayC}.
\bigskip

Assume $\lambda,\mu$ belong to $\mathcal{P}^{\diamondsuit}$.\ We write
$\mu\subset\lambda$ when the Young diagram of $\mu$ is contained in the one of
$\lambda$. In that case, the skew Young diagram $\lambda/\mu$ is obtained from
$\lambda$ by deleting the boxes appearing in $\mu$.

\subsection{Tensor powers of the natural representation}

Each algebra $\mathfrak{g}\mathfrak{=}\mathfrak{gl}(n),\mathfrak{gl}(m,n)$ and
$\mathfrak{q}(n)$ can be realized as a matrix algebra.\ They thus admit a
natural representation $V^{\diamondsuit}$ which is the vector representation
of the underlying matrix algebra.\ For any $\ell\geq0,$ the tensor power
$(V^{\diamondsuit})^{\otimes\ell}$ is a semisimple representation for
$\mathfrak{g}$.\ This means that $(V^{\diamondsuit})^{\otimes\ell}$ decomposes
into a direct sum of irreducible representations
\begin{equation}
(V^{\diamondsuit})^{\otimes\ell}\simeq\bigoplus_{\lambda\in P_{+}
}V^{\diamondsuit}(\lambda)^{\oplus f_{\lambda}^{\diamondsuit}}
\end{equation}
where for any $\lambda\in P_{+},$ the module $V^{\diamondsuit}(\lambda)$ is
the irreducible module with highest weight $\lambda$ and multiplicity
$f_{\lambda}^{\diamondsuit}$ in $(V^{\diamondsuit})^{\otimes\ell}$. In
general, we cannot realize all the irreducible highest weight modules as
irreducible components in a tensor product $(V^{\diamondsuit})^{\otimes\ell}.$
More precisely, we have the following proposition.

\begin{proposition}
For any $\lambda\in P_{+},$ the module $V^{\diamondsuit}(\lambda)$ appears as
an irreducible component in a tensor product $(V^{\diamondsuit})^{\otimes\ell
}$ if and only if $\lambda\in\mathcal{P}^{\diamondsuit}$ and $\left|
\lambda\right|  =\ell$.
\end{proposition}

When $\mu\in\mathcal{P}^{\diamondsuit}$, we also define the multiplicities
$f_{\lambda/\mu}^{\diamondsuit}$ by
\begin{equation}
V^{\diamondsuit}(\mu)\otimes(V^{\diamondsuit})^{\otimes\ell}\simeq
\bigoplus_{\lambda\in P_{+}}V^{\diamondsuit}(\lambda)^{\oplus f_{\lambda/\mu
}^{\diamondsuit}}. \label{skew_multi}
\end{equation}
Set $\ell^{\prime}=\left\vert \mu\right\vert $. Since $V^{\diamondsuit}(\mu)$
appears as an irreducible component of $(V^{\diamondsuit})^{\otimes
\ell^{\prime}},$ one gets $f_{\lambda/\mu}^{\diamondsuit}\neq0$ if and only if
$\lambda$ appears as an irreducible component of $(V^{\diamondsuit}
)^{\otimes\ell+\ell^{\prime}}$. In this situation, we have $\lambda
\in\mathcal{P}^{\diamondsuit}$ and $\left\vert \lambda\right\vert =\ell
+\ell^{\prime}$.\ When $\ell=1$ we have
\begin{equation}
V^{\diamondsuit}(\mu)\otimes V^{\diamondsuit}\simeq\bigoplus_{\mu
\rightsquigarrow\lambda}V^{\diamondsuit}(\lambda) \label{ten_by_V}
\end{equation}
where $\mu\rightsquigarrow\lambda$ means that the sum is over all the
partitions $\lambda\in\mathcal{P}^{\diamondsuit}$ obtained from $\mu
\in\mathcal{P}^{\diamondsuit}$ by adding one box. More generally if $\kappa
\in\mathcal{P}^{\diamondsuit},$ we set
\begin{equation}
V^{\diamondsuit}(\mu)\otimes V^{\diamondsuit}(\kappa)=\bigoplus_{\lambda
\in\mathcal{P}^{\diamondsuit}}V^{\diamondsuit}(\lambda)^{\oplus m_{\kappa,\mu
}^{\lambda,\diamondsuit}}. \label{LRcoef}
\end{equation}
Recall that $V^{\diamondsuit}(\mu)\otimes V^{\diamondsuit}(\kappa)$ is
isomorphic to $V^{\diamondsuit}(\kappa)\otimes V^{\diamondsuit}(\mu)$ and
therefore $m_{\kappa,\mu}^{\lambda,\diamondsuit}=m_{\mu,\kappa}^{\lambda
,\diamondsuit}$. Observe that $m_{\mu,(1)}^{\lambda,\diamondsuit}=1$ if and only
if $\mu\rightsquigarrow\lambda.$

\section{Combinatorics of tableaux}

\label{sec_ComTab}In this section, we review the different notions of tableaux
which are relevant for the representation theory of $\mathfrak{gl}
(n),\mathfrak{gl}(m,n)$ and $\mathfrak{q}(n).$ We also recall the associated
insertion schemes which are essential to define the generalized Pitman
transform in an elementary way. These notions are very classical for
$\mathfrak{gl}(n)$ (see \cite{Fu}); for $\mathfrak{gl}(m,n),$ they were
introduced by Benkart, Kang and Kashiwara in \cite{BKK}; for $\mathfrak{q}(n)
$ this has been developed very recently in \cite{GHSKM}.

\subsection{Characters and tableaux}

\subsubsection{Semistandard $\mathfrak{gl}(n)$-tableaux}

Consider $\lambda\in\mathcal{P}^{\emptyset}$.\ A (semistandard) $\mathfrak{gl}
(n)$-tableau of shape $\lambda$ is a filling (let us call it $T $) of the
Young diagram associated to $\lambda$ by letters of the ordered alphabet
$\mathcal{A}_{n}=\{1<2<\cdots<n\}$ such that the rows of $T$ weakly increase
from left to right and its columns strictly increase from top to bottom (see
Example \ref{ex_A}). We denote by $T^{\emptyset}(\lambda)$ the set of all
$\mathfrak{gl}(n)$-tableaux of shape $\lambda.\;$We define the reading of
$T\in T^{\emptyset}(\lambda)$ as the word $\mathrm{w}(T)$ of $\mathcal{A}
_{n}^{\ast}$ obtained by reading the rows of $T$ from right to left and then
top to bottom.

The weight of a word $w\in\mathcal{A}_{n}^{\ast}$ is the $n$-tuple
$\mathrm{wt}(w)=(\mu_{1},\ldots,\mu_{n})$ where for any $i=1,\ldots, n $ the
nonnegative integer $\mu_{i}$ is the number of letters $i$ in $w$ for any
$i=1, \cdots, n.\;$The weight $\mathrm{wt}(T)$ of $T\in T^{\emptyset}
(\lambda)$ is then defined as the weight of its reading $\mathrm{w}(T)$. The
Schur function $s_{\lambda}^{\emptyset}$ is the character of $V^{\emptyset
}(\lambda)$. This is a symmetric polynomial in the variables $x_{1}
,\ldots,x_{n}$ which can be expressed as a generating series over
$T^{\emptyset}(\lambda),$ namely we have
\begin{equation}\label{schur-funct-gn}
s_{\lambda}^{\emptyset}(x)=\sum_{T\in T^{\emptyset}(\lambda)}x^{\mathrm{wt}
(T)}.
\end{equation}
According to the Weyl character formula, we have
\begin{equation}
s_{\lambda}^{\emptyset}(x)=\frac{1}{\prod_{1\leq i<j\leq n}(1-\frac{x_{j}
}{x_{i}})}\sum_{\sigma\in S_{n}}\varepsilon(\sigma)x^{\sigma(\lambda
+\rho)-\rho}=\frac{1}{\prod_{1\leq i<j\leq n}(x_{i}-x_{j})}\sum_{\sigma\in
S_{n}}\varepsilon(\sigma)x^{\sigma(\lambda+\rho)} \label{Weyl}
\end{equation}
where $S_{n}$ is the symmetric group of rank $n,$ $\varepsilon(\sigma)$ the
signature of $\sigma,$ $\rho=(n-1,n-2,\ldots,0)\in\mathbb{Z}^{n}$ and $S_{n}$
acts on $\mathbb{Z}^{n}$ by permutation of the coordinates.

\subsubsection{Semistandard $\mathfrak{gl}(m,n)$-tableaux}

Consider $\lambda\in\mathcal{P}^{h}$.\ A (semistandard) $\mathfrak{gl}
(m,n)$-tableau of shape $\lambda$ is a filling $T$ of the hook Young diagram
associated to $\lambda$ by letters of the ordered alphabet
\[
\mathcal{A}_{m,n}=\{\overline{m}<\overline{m-1}<\cdots<\overline{1}
<1<2<\cdots<n\}
\]
such that the rows of $T$ weakly increase from left to right with no
repetition of unbarred letters permitted and its columns increase from top to
bottom with no repetition of barred letters permitted (see Example
\ref{ex_h}). We denote by $T^{h}(\lambda)$ the set of all $\mathfrak{gl}
(m,n)$-tableaux of shape $\lambda.\;$We define the reading of $T\in
T^{h}(\lambda)$ as the word $\mathrm{w}(T)$ of $\mathcal{A}_{m,n}^{\ast}$
obtained by reading the rows of $T$ from right to left and top to bottom.

The weight of a word $w\in\mathcal{A}_{m,n}$ is the $(m+n)$-tuple
$\mathrm{wt}(w)=(\mu_{\overline{m}},\ldots,\mu_{\overline{1}}\mid\nu
_{1},\ldots,\nu_{n})$ where for any $i=1,\ldots,m $ and $j=1,\ldots,n,$ the
nonnegative integer $\mu_{\overline{i}}$ is the number of letters
$\overline{i}$ in $w$ and $\nu_{j}$ is the number of letters $j$ in
$w$. The weight of $T\in T^{h}(\lambda)$ is then defined as the weight of its
reading $\mathrm{w}(T)$. The Schur function $s_{\lambda}^{h}$ is the character of the
irreducible representation $V^{h}(\lambda)$ of $\mathfrak{gl}(m,n)$. This is a
polynomial in the variables $x_{\overline{m}},\ldots,x_{x_{\overline{1}}
},x_{1},\ldots,x_{n}$. It admits a nice expression as a generating series
over $T^{h}(\lambda)$, namely we have
\begin{equation}\label{schur-funct-gmn}
s_{\lambda}^{h}(x)=\sum_{T\in T^{h}(\lambda)}x^{\mathrm{wt}(T)}.
\end{equation}
For a general highest weight $\mathfrak{gl}(m,n)$-module there is no simple
Weyl character formula. Nevertheless for the irreducible modules
$V^{h}(\lambda)$ with $\lambda\in^{h}$, such a formula exists due to Berele,
Regev, and Serge'ev (see \cite{Hu}). Consider $\lambda\in\mathcal{P}^{h}:$ the
character $s_{\lambda}^{h}$ of $V^{h}(\lambda)$ is given by
\begin{equation}
s_{\lambda}^{h}(x)=\frac{\prod_{(i,j)\in\lambda}(1+\frac{x_{j}}{x_{\overline
{i}}})}{\prod_{\overline{m}\leq\overline{i}<\overline{j}\leq\overline{1}
}(1-\frac{x_{\overline{j}}}{x_{\overline{i}}})\prod_{1\leq r<s\leq n}
(1-\frac{x_{s}}{x_{r}})}\sum_{w\in S_{m}\times S_{n}}\varepsilon
(w)\,x^{w(\lambda+\rho_{+})-\rho_{+}} \label{superWeyl}
\end{equation}
where $\rho_{+}=(m-1,\ldots,1,0\mid n-1,\ldots,1,0)$ and $(i,j)\in\lambda$
means that the hook Young diagram associated to $\lambda$ has a box at the
intersection of its row $i\in\{1,\ldots,m\}$ and its column $j\in
\{1,\ldots,n\}.$

\subsubsection{Semistandard decomposition $\mathfrak{q}(n)$-tableaux}

Let us first give a definition. We say that a nonempty word $w=x_{1}\cdots x_{\ell}\in\mathcal{A}_{n}^{\ast}$
is a \textit{hook word} if there exists $1\leq k\leq\ell$ such that $x_{1}\geq
x_{2}\geq\cdots\geq x_{k}<x_{k+1}<\cdots<x_{\ell}$. Each hook word can be decomposed $w=w_{\downarrow}w_{\uparrow}$ where by convention,
the \emph{decreasing part} $w_{\downarrow}=x_{1}\geq x_{2}\geq\cdots\geq x_{k}$ is nonempty. The \emph{increasing part}
$w_{\uparrow}=x_{k+1}<\cdots<x_{\ell}$ is possibly empty. In particular when
$w=x_{1}\cdots x_{\ell}$ is such that $x_{1}<\cdots<x_{\ell}$, we have
$w_{\downarrow}=x_{1}$ and $w_{\uparrow}=x_{2}\cdots x_{\ell}$.

Consider a strict partition $\lambda\in\mathcal{P}^{s}$.\ A (semistandard)
$\mathfrak{q}(n)$-tableau of shape $\lambda$ is a filling $T$ of
the shifted Young diagram associated to $\lambda$ by letters of $\mathcal{A}
_{n}=\{1<2<\cdots<n\}$ such that for $i=1,\ldots,n,$

\begin{enumerate}
\item the word $w_{i}$ formed by reading the $i$-th row of $T$ \emph{from left
to right} is a hook word (of length $\lambda_{i}$),

\item $w_{i}$ is a hook subword of maximal length in $w_{i+1}w_{i}$ (see
Example \ref{ex_s}).
\end{enumerate}

We denote by $T^{s}(\lambda)$ the set of all $\mathfrak{q}(n)$-tableaux of
shape $\lambda.\;$The reading of $T\in T^{s}(\lambda)$ is the word
$\mathrm{w}(T)=w_{n}\cdots w_{2}w_{1}$ of $\mathcal{A}_{n}^{\ast}.$ We define
the weight of $T$ as the weight of its reading. The Schur function
$s_{\lambda}^{s}$ is defined as the generating series
\begin{equation}\label{schur-funct-qn}
s_{\lambda}^{s}(x)=\sum_{T\in T^{s}(\lambda)}x^{\mathrm{wt}(T)}.
\end{equation}
This is not the original combinatorial definition of the Schur function
which was given in terms of different tableaux called \emph{shifted Young tableaux}.
Nevertheless, according to Theorem 2.17 in \cite{Sera} and Remark 2.6 in
\cite{GHSKM} there exists a weight preserving bijection between the set of
shifted Young tableaux of shape $\lambda$ and the set of $\mathfrak{q}
(n)$-tableaux with the same shape. Set $d(\lambda)$ for the depth of $\lambda
$, that is for the number of nonzero coordinates in $\lambda$. The Schur
function admits a Weyl type expression, namely we have.
\[
s_{\lambda}^{s}(x)
=\sum_{\sigma\in
S_{n}/S_{\lambda}}\sigma\left(  x^{\lambda}\prod_{1\leq i\leq d(\lambda)}
\prod_{i<j\leq n}\left(  \dfrac{x_{i}+x_{j}}{x_{i}-x_{j}}\right)  \right)
\]
where $S_{\lambda}$ is the stabilizer of $\lambda$ under the action of $S_{n}
$. Thus $S_{\lambda}$ is either isomorphic to $S_{n-d(\lambda)}$ when $d(\lambda
)<n-1$ or it reduces to $\{id\}$ otherwise.

\noindent Assume $\lambda\in\mathcal{P}^{s}$ is such that $d(\lambda)=n.$ In
that case we have $S_{\lambda}=\{id\}$ and
\[
s_{\lambda}^{s}(x)=\left(  \sum_{\sigma\in S_{n}}\varepsilon(\sigma
)x^{\sigma(\lambda)}\right)  \times\prod_{1\leq i<j\leq n}\left(  \dfrac
{x_{i}+x_{j}}{x_{i}-x_{j}}\right)
\]
since the function $\displaystyle \prod_{1\leq i<j\leq n}x_{i}+x_{j}$ is symmetric and 
$\displaystyle \prod_{1\leq i<j\leq n}
x_{i}-x_{j}$  is antisymmetric. Then $\lambda
-\rho\in\mathcal{P}^{\emptyset}$. It follows from (\ref{Weyl}) and the above
equality that we have the identity $^{(}$\footnote{we thank here Marc Van
Leeuwen who pointed out to us this property.}$^{)}$
\begin{equation}
s_{\lambda}^{s}(x)=s_{\lambda-\rho}^{\emptyset}(x)\prod_{1\leq i<j\leq
n}\left(  x_{i}+x_{j}\right)  . \label{rela}
\end{equation}

\subsection{Insertion schemes}

\label{sub_sec_insert}To make our notation consistent, we set $\mathcal{A}
^{\emptyset}=\mathcal{A}^{s}=\mathcal{A}_{n}$ and $\mathcal{A}^{h}
=\mathcal{A}_{m,n}.$

\subsubsection{On $\mathfrak{gl}(n)$-tableaux}

Let $T$ be a $\mathfrak{gl}(n)$-tableau of shape $\lambda\in\mathcal{P}
^{\emptyset}.\;$We write $T=C_{1}\cdots C_{s}$ as the juxtaposition of its
columns. Consider $x\in\mathcal{A}_{n}.\;$We denote by $x\rightarrow T$ the
tableau obtained by applying the following recursive procedure:

\begin{enumerate}
\item Assume $T=\emptyset,$ then $x\rightarrow T$ is the tableau with one box
filled by $x$.

\item Assume $C_{1}$ is nonempty.

\begin{enumerate}
\item If all the letters of $C_{1}$ are less than $x,$ the tableau
$x\rightarrow T$ is obtained from $T$ by adding one box filled by $x$ at the
bottom of $C_{1}.$

\item Otherwise, let $y=\min\{t\in C_{1}\mid t\geq x\}.\;$Write $C_{1}
^{\prime}$ for the column obtained by replacing $y$ by $x$ in $C_{1}.$ Then
$x\rightarrow T=C_{1}^{\prime}(y\rightarrow C_{2}\cdots C_{s})$ is defined as
the juxtaposition of $C_{1}^{\prime}$ with the tableau obtained by inserting
$y$ in the remaining columns.
\end{enumerate}
\end{enumerate}

One easily verifies that in any case $x\rightarrow T$ is a $\mathfrak{gl}
(n)$-tableau. More generally, for any word $w=x_{1}x_{2}\cdots x_{\ell}
\in\mathcal{A}_{n}^{\ast}$, we define the $\mathfrak{gl}(n)$-tableau
$P^{\emptyset}(w)$ setting
\begin{equation}
P^{\emptyset}(w)=x_{\ell}\rightarrow(\cdots(x_{2}\rightarrow(x_{1}
\rightarrow\emptyset))). \label{defP(w)}
\end{equation}

\begin{example}
\label{ex_A}With $n\geq4$ and $w=232143,$ we obtain the following sequences of
tableaux:
\[
\begin{tabular}
[c]{|l|}\hline
$\mathtt{2}$\\\hline
\end{tabular}
\ \text{, }
\begin{tabular}
[c]{|l|}\hline
$\mathtt{2}$\\\hline
$\mathtt{3}$\\\hline
\end{tabular}
\ \text{, }
\begin{tabular}
[c]{|l|l}\hline
$\mathtt{2}$ & \multicolumn{1}{|l|}{$\mathtt{2}$}\\\hline
$\mathtt{3}$ & \\\cline{1-1}
\end{tabular}
\ \text{,
\begin{tabular}
[c]{|l|ll}\hline
$\mathtt{1}$ & $\mathtt{2}$ & \multicolumn{1}{|l|}{$\mathtt{2}$}\\\hline
$\mathtt{3}$ &  & \\\cline{1-1}
\end{tabular}
, }
\begin{tabular}
[c]{|l|ll}\hline
$\mathtt{1}$ & $\mathtt{2}$ & \multicolumn{1}{|l|}{$\mathtt{2}$}\\\hline
$\mathtt{3}$ &  & \\\cline{1-1}
$\mathtt{4}$ &  & \\\cline{1-1}
\end{tabular}
\ \text{, }
\begin{tabular}
[c]{|l|ll}\hline
$\mathtt{1}$ & $\mathtt{2}$ & \multicolumn{1}{|l|}{$\mathtt{2}$}\\\hline
$\mathtt{3}$ & $\mathtt{3}$ & \multicolumn{1}{|l}{}\\\cline{1-2}
$\mathtt{4}$ &  & \\\cline{1-1}
\end{tabular}
\ =P^{\emptyset}(w).
\]

\end{example}

\subsubsection{On $\mathfrak{gl}(m,n)$-tableaux}

Let $T=C_{1}\cdots C_{s}$ be a $\mathfrak{gl}(m,n)$-tableau of shape
$\lambda\in\mathcal{P}^{h}.\;$Consider $x\in\mathcal{A}_{m,n}.\;$We denote by
$x\rightarrow T$ the tableau obtained by applying the following procedure:

\begin{enumerate}
\item Assume $T=\emptyset,$ then $x\rightarrow T$ is the tableau with one box
filled by $x$.

\item Assume $C_{1}$ is nonempty and $x$ is a barred letter.

\begin{enumerate}
\item If all the letters of $C_{1}$ are less than $x,$ the tableau
$x\rightarrow T$ is obtained from $T$ by adding one box filled by $x$ at the
bottom of $C_{1}.$

\item Otherwise, let $y=\min\{t\in C_{1}\mid t\geq x\}.\;$Write $C_{1}
^{\prime}$ for the column obtained by replacing $y$ by $x$ in $C_{1}$. If $y$
appears at least twice in $C_{1}$ (this can happen when $y$ is unbarred), we
replace by $x$ the letter $y$ appearing in the highest position in $C_{1}
$.\ Then $x\rightarrow T=C_{1}^{\prime}(y\rightarrow C_{2}\cdots C_{s}).$
\end{enumerate}

\item Assume $C_{1}$ is nonempty and $x$ is an unbarred letter.

\begin{enumerate}
\item If all the letters of $C_{1}$ are less or equal to $x,$ the tableau
$x\rightarrow T$ is obtained from $T$ by adding one box filled by $x$ at the
bottom of $C_{1}.$

\item Otherwise, let $y=\min\{t\in C_{1}\mid t>x\}.\;$Write $C_{1}^{\prime}$
for the column obtained by replacing $y$ by $x$ in $C_{1}$. Similarly, if $y$
appears at least twice in $C_{1}$ we replace by $x$ the letter $y$ appearing
in the highest position in $C_{1}$.\ Then $x\rightarrow T=C_{1}^{\prime
}(y\rightarrow C_{2}\cdots C_{s}).$
\end{enumerate}
\end{enumerate}

One verifies that in any cases $x\rightarrow T$ is a $\mathfrak{gl}
(m,n)$-tableau. For any word $w=x_{1}x_{2}\cdots x_{\ell}\in\mathcal{A}
_{m,n}^{\ast}$, we define recursively the $\mathfrak{gl}(m,n)$-tableau
$P^{h}(w)$ as in (\ref{defP(w)}).

\begin{example}
\label{ex_h}With $(m,n)=(2,3)$ and $w=\bar{2}3\bar{2}\bar{1}32\bar{1}2,$ we
obtain the following sequence of tableaux:
\[
\begin{tabular}
[c]{|l|}\hline
$\mathtt{\bar{2}}$\\\hline
\end{tabular}
\text{, }
\begin{tabular}
[c]{|l|}\hline
$\mathtt{\bar{2}}$\\\hline
$\mathtt{3}$\\\hline
\end{tabular}
\text{, }
\begin{tabular}
[c]{|l|l}\hline
$\mathtt{\bar{2}}$ & \multicolumn{1}{|l|}{$\mathtt{\bar{2}}$}\\\hline
$\mathtt{3}$ & \\\cline{1-1}
\end{tabular}
\text{, }
\begin{tabular}
[c]{|l|l|}\hline
$\mathtt{\bar{2}}$ & $\mathtt{\bar{2}}$\\\hline
$\mathtt{\bar{1}}$ & $\mathtt{3}$\\\hline
\end{tabular}
\text{, }
\begin{tabular}
[c]{|l|l}\hline
$\mathtt{\bar{2}}$ & \multicolumn{1}{|l|}{$\mathtt{\bar{2}}$}\\\hline
$\mathtt{\bar{1}}$ & \multicolumn{1}{|l|}{$\mathtt{3}$}\\\hline
$\mathtt{3}$ & \\\cline{1-1}
\end{tabular}
\text{, }
\begin{tabular}
[c]{|l|l|}\hline
$\mathtt{\bar{2}}$ & $\mathtt{\bar{2}}$\\\hline
$\mathtt{\bar{1}}$ & $\mathtt{3}$\\\hline
$\mathtt{2}$ & $\mathtt{3}$\\\hline
\end{tabular}
\text{, }
\begin{tabular}
[c]{|l|l|l}\hline
$\mathtt{\bar{2}}$ & $\mathtt{\bar{2}}$ & \multicolumn{1}{|l|}{$\mathtt{3}$
}\\\hline
$\mathtt{\bar{1}}$ & $\mathtt{\bar{1}}$ & \\\cline{1-2}\cline{2-2}
$\mathtt{2}$ & $\mathtt{3}$ & \\\cline{1-2}
\end{tabular}
\text{, }
\begin{tabular}
[c]{|l|ll}\hline
$\mathtt{\bar{2}}$ & $\mathtt{\bar{2}}$ & \multicolumn{1}{|l|}{$\mathtt{3}$
}\\\hline
$\mathtt{\bar{1}}$ & $\mathtt{\bar{1}}$ & \multicolumn{1}{|l}{}\\\cline{1-2}
\cline{2-2}
$\mathtt{2}$ & $\mathtt{3}$ & \multicolumn{1}{|l}{}\\\cline{1-2}
$\mathtt{2}$ &  & \\\cline{1-1}
\end{tabular}
=P^{h}(w).
\]

\end{example}

\subsubsection{On $\mathfrak{q}(n)$-tableaux}

Let $T=L_{1}\cdots L_{k}$ be a $\mathfrak{q}(n)$-tableau of shape $\lambda
\in\mathcal{P}^{s}.\;$Here we regard $T$ as the juxtaposition of its rows
(written by decreasing lengths) rather as the juxtaposition of its
columns.\ Consider $x\in\mathcal{A}_{n}.\;$We denote by $x\rightarrow T$ the
tableau obtained by applying the following procedure (which implies in
particular that, at each step and for any row $L_{i}$ the word $w(L_{i})$ is
a hook word):

\begin{enumerate}
\item Assume $T=\emptyset,$ then $x\rightarrow T$ is the tableau with one box
filled by $x$.

\item Assume $L_{1}$ is nonempty and write $w=w_{\downarrow}w_{\uparrow}$ the
decomposition of $\mathrm{w}(L_{1})$ as decreasing and increasing subwords.

\begin{enumerate}
\item If $wx$ is a hook word, then $x\rightarrow T$ is the tableau obtained
from $T$ by adding one box filled by $x$ at the right end of $L_{1}.$

\item Otherwise, $w_{\uparrow}\neq\emptyset$ and $y=\min\{t\in w_{\uparrow
}\mid t\geq x\}$ exists. We first replace $y$ by $x$ in $w_{\uparrow}.$ Now
let $z=\max\{t\in w_{\downarrow}\mid t<y\}.\;$We replace $z$ by $y$ in
$w_{\downarrow}.$ Write $L_{1}^{\prime}$ for the row so obtained.\ Then
$x\rightarrow T=L_{1}^{\prime}(z\rightarrow L_{2}\cdots L_{k}).$
\end{enumerate}
\end{enumerate}

One also verifies that this gives a $\mathfrak{q}(n)$-tableau. For any word
$w=x_{1}x_{2}\cdots x_{\ell}\in\mathcal{A}_{n}^{\ast}$, we define recursively
the $\mathfrak{q}(n)$-tableau $P^{s}(w)$ as in (\ref{defP(w)}).

\begin{example}
\label{ex_s}With $n=4$ and $w=232145331,$ we obtain the following sequence of
tableaux:
\begin{gather*}
\begin{tabular}
[c]{|l|}\hline
$2$\\\hline
\end{tabular}
\ \text{, }
\begin{tabular}
[c]{|l|l|}\hline
$2$ & $\mathbf{3}$\\\hline
\end{tabular}
\ \text{, }
\begin{tabular}
[c]{c|c|}\hline
\multicolumn{1}{|c|}{$3$} & $2$\\\hline
& $2$\\\cline{2-2}
\end{tabular}
\ \text{,
\begin{tabular}
[c]{c|c|c}\hline
\multicolumn{1}{|c|}{$3$} & $2$ & \multicolumn{1}{|c|}{$1$}\\\hline
& $2$ & \\\cline{2-2}
\end{tabular}
,
\begin{tabular}
[c]{c|c|cc}\hline
\multicolumn{1}{|c|}{$3$} & $2$ & $1$ & \multicolumn{1}{|c|}{$\mathbf{4}$
}\\\hline
& $2$ &  & \\\cline{2-2}
\end{tabular}
, }
\begin{tabular}
[c]{c|c|ccc}\hline
\multicolumn{1}{|c|}{$3$} & $2$ & $1$ & \multicolumn{1}{|c}{$\mathbf{4}$} &
\multicolumn{1}{|c|}{$\mathbf{5}$}\\\hline
& $2$ &  &  & \\\cline{2-2}
\end{tabular}
\ \text{, }\\
\begin{tabular}
[c]{c|c|c|cc}\hline
\multicolumn{1}{|c|}{$4$} & $2$ & $1$ & $\mathbf{3}$ &
\multicolumn{1}{|c|}{$\mathbf{5}$}\\\hline
& $2$ & $\mathbf{3}$ &  & \\\cline{2-3}
\end{tabular}
\ \text{, }
\begin{tabular}
[c]{cc|c|cc}\hline
\multicolumn{1}{|c}{$4$} & \multicolumn{1}{|c|}{$3$} & $1$ & $\mathbf{3}$ &
\multicolumn{1}{|c|}{$\mathbf{5}$}\\\hline
& \multicolumn{1}{|c|}{$3$} & $2$ &  & \\\cline{2-3}\cline{3-3}
&  & $2$ &  & \\\cline{3-3}
\end{tabular}
\ \text{, }
\begin{tabular}
[c]{cc|c|cc}\hline
\multicolumn{1}{|c}{$4$} & \multicolumn{1}{|c|}{$3$} & $3$ & $1$ &
\multicolumn{1}{|c|}{$\mathbf{5}$}\\\hline
& \multicolumn{1}{|c|}{$3$} & $2$ & $1$ & \multicolumn{1}{|c}{}\\\cline{2-4}
\cline{3-3}
&  & $2$ &  & \\\cline{3-3}
\end{tabular}
\ =P^{s}(w).
\end{gather*}
where we have indicated in bold the increasing part of each row.
\end{example}

\subsection{Robinson-Schensted correspondence}

\label{subsec_RSK} For any word $w=x_{1}\cdots x_{\ell}\in(\mathcal{A}
^{\diamondsuit})^{\ell}$ and any $k=1,\ldots\ell,$ let $\lambda^{(k)} $ be the
shape of the tableau $P^{\diamondsuit}(x_{1}\cdots x_{k})$; the shape
$\lambda^{(k)}$ is obtained from $\lambda^{(k-1)}$ by adding one box we denote
by $b_{k}.\;$ The recording tableau $Q^{\diamondsuit}(w)$ of shape
$\lambda^{(\ell)}$ is obtained by filling each box $b_{k}$ with the letter
$k$. Observe that $Q^{\diamondsuit}(w)$ is a standard tableau: it contains
exactly once all the integers $1,\ldots,\ell$, its rows strictly increase from
left to right and its columns strictly increase from top to bottom. Note also
that the datum of a standard tableau with $\ell$ boxes is equivalent to that
of a sequence of shapes $(\lambda^{(1)},\ldots,\lambda^{(\ell)})\in
(\mathcal{P}^{\diamondsuit})^{\ell}$ such that $\lambda^{(1)}=(1)$ and for any
$k=1,\ldots,\ell$, the shape $\lambda^{(k )}$ is obtained by adding one box to
$\lambda^{(k-1)}$.

\begin{examples}
From the previous examples, we derive

$Q^{\emptyset}(232143)=
\begin{tabular}
[c]{|l|ll}\hline
$\mathtt{1}$ & $\mathtt{3}$ & \multicolumn{1}{|l|}{$\mathtt{4}$}\\\hline
$\mathtt{2}$ & $\mathtt{6}$ & \multicolumn{1}{|l}{}\\\cline{1-2}
$\mathtt{5}$ &  & \\\cline{1-1}
\end{tabular}
,$ $Q^{h}(\bar{2}3\bar{2}\bar{1}32\bar{1}2)=
\begin{tabular}
[c]{|l|ll}\hline
$\mathtt{1}$ & $\mathtt{3}$ & \multicolumn{1}{|l|}{$\mathtt{7}$}\\\hline
$\mathtt{2}$ & $\mathtt{4}$ & \multicolumn{1}{|l}{}\\\cline{1-2}\cline{2-2}
$\mathtt{5}$ & $\mathtt{6}$ & \multicolumn{1}{|l}{}\\\cline{1-2}
$\mathtt{8}$ &  & \\\cline{1-1}
\end{tabular}
\ $and $Q^{s}(23214433)=
\begin{tabular}
[c]{cc|c|cc}\hline
\multicolumn{1}{|c}{$1$} & \multicolumn{1}{|c|}{$2$} & $4$ & $5$ &
\multicolumn{1}{|c|}{$6$}\\\hline
& \multicolumn{1}{|c|}{$3$} & $7$ & $9$ & \multicolumn{1}{|c}{}\\\cline{2-4}
\cline{3-3}
&  & $8$ &  & \\\cline{3-3}
\end{tabular}
\ .$
\end{examples}

By a $\diamondsuit$-tableau $\diamondsuit=\emptyset,h,s$, we mean a tableau
for $\mathfrak{gl}(n),\mathfrak{gl}(m,n),\mathfrak{q}(n)$.\ We can now state
the Robinson-Schensted correspondence for $\mathfrak{gl}(n)$ (see \cite{Fu})
and its generalizations for $\mathfrak{gl}(m,n)$ and $\mathfrak{q}(n)$
obtained in \cite{BKK} and \cite{GHSKM}, respectively. For any $\ell
\geq0$, write $\mathcal{U}_{\ell}^{\diamondsuit}$
($\diamondsuit\in\{\emptyset,h,s\}$) for the set of pairs $(P,Q)$ where $P$ is
a $\diamondsuit$-tableau and $Q$ a standard tableau \emph{with the same shape}
as $P$ containing $\ell$ boxes.

Consider $\lambda\in\mathcal{P}^{\diamondsuit}$ and assume $\left\vert
\lambda\right\vert =\ell.$ Given $T$ a standard tableau of shape $\lambda, $
we set
\[
B^{\diamondsuit}(T)=\{w\in(\mathcal{A}^{\diamondsuit})^{\ell}\mid
Q^{\diamondsuit}(w)=T\}.
\]

One may state the

\begin{theorem}
\ \label{Th_RSK} \cite{GHSKM} Fix $\diamondsuit\in\{\emptyset,h,s\}$.

\begin{enumerate}
\item The map $\displaystyle
\left\{  \theta_{\ell}^{\diamondsuit}:
\begin{array}
[c]{l}
(\mathcal{A}^{\diamondsuit})^{\ell}\rightarrow\mathcal{U}_{\ell}
^{\diamondsuit}\\
w\mapsto(P^{\diamondsuit}(w),Q^{\diamondsuit}(w))
\end{array}
\right.  $ is a one-to-one correspondence.

In particular, the map $P^{\diamondsuit}$ restricts to a weight preserving
bijection $P^{\diamondsuit}:B^{\diamondsuit}(T)\longleftrightarrow
T^{\diamondsuit}(\lambda).$

\item For any $\lambda\in\mathcal{P}^{\diamondsuit}$, the multiplicity
$f_{\lambda}^{\diamondsuit}$ is equal to the number of standard $\diamondsuit
$-tableaux of shape~$\lambda.$
\end{enumerate}
\end{theorem}

Given $\lambda,\mu$ in $\mathcal{P}^{\diamondsuit}$ regarded as Young diagrams
such that $\mu_{i}\leq\lambda_{i}$ for any $i$ with $\mu_{i}>0$, we denote by
$\lambda/\mu$ the skew Young diagram obtained by deleting in $\lambda$ the
boxes of $\mu.\;$By a standard tableau of shape $\lambda/\mu$ with $\ell$
boxes, we mean a filling of $\lambda/\mu$ by the letters of $\{1,\ldots
,\ell\}$ whose rows and columns strictly increase from left to right and top
to bottom, respectively. By a skew $\diamondsuit$-tableau of shape
$\lambda/\mu,$ we mean a filling of $\lambda/\mu$ by letters of $\mathcal{A}
^{\diamondsuit}$ whose rows and columns satisfy the same conditions as for the
ordinary $\diamondsuit$-tableaux. The following Proposition will follow from
the Littelwood-Richardson rules proved in \cite{Fu} for $\mathfrak{gl}(n)$, in \cite{KK} for $\mathfrak{gl}(m,n)$
and in \cite{GHSKM} for $\mathfrak{q}(n)$.\ We postpone its proof to the Appendix.

\begin{proposition}
\label{Cor_RSK} \ 

\begin{enumerate}
\item Given $\lambda,\mu$ in $\mathcal{P}^{\diamondsuit}$ such that
$\mu\subset\lambda,$ the multiplicity $f_{\lambda/\mu}^{\diamondsuit}$ defined
in (\ref{skew_multi}) is equal to the number of standard tableaux of shape
$\lambda/\mu$.

\item Given $\lambda,\kappa,\mu$ in $\mathcal{P}^{\diamondsuit}$, we have
$m_{\kappa,\mu}^{\lambda,\diamondsuit}\leq K_{\mu,\lambda-\kappa
}^{\diamondsuit}$ where $K_{\mu,\lambda-\kappa}^{\diamondsuit}$ is the weight
multiplicity defined in (\ref{Kostka}) and $m_{\kappa,\mu}^{\lambda
,\diamondsuit}$ the tensor multiplicity defined in (\ref{LRcoef}).
\end{enumerate}
\end{proposition}

\noindent\textbf{Remarks:} $\mathrm{(i)}$ It follows from the definition of Schur functions (\ref{schur-funct-gn}), (\ref{schur-funct-gmn}) and (\ref{schur-funct-qn}), and from 1 of Theorem
\ref{Th_RSK} that for $T$ a standard tableau of shape $\lambda\in
\mathcal{P}^{\diamondsuit},$ we have
\[
\sum_{w\in B^{\diamondsuit}(T)}x^{\mathrm{wt}(w)}=s_{\lambda}^{\diamondsuit
}(x).
\]

\noindent$\mathrm{(ii)}$ Using Kashiwara crystal basis theory, we can
obtain a stronger version of the previous theorem.\ For any $\diamondsuit
\in\{\emptyset,h,s\}$, the set $B^{\diamondsuit}(T)$ has a crystal structure :
namely, it can be endowed with the structure of an oriented graph (depending
on $\diamondsuit$) with arrows colored by integers.\ Such a structure can also
be defined on the set of words $(\mathcal{A}^{\diamondsuit})^{\ell}$; one
then shows that the sets $T^{\diamondsuit}(T)$ where $T$ runs over the set of
standard tableaux with $\ell$ boxes are the connected components of
$(\mathcal{A}^{\diamondsuit})^{\ell}.$ Moreover the bijection of assertion 1.
of the theorem is a graph isomorphism, that is compatible with this crystal
structure. For $\mathfrak{gl}(n)$ the crystal basis theory is now a classical
tool in representation theory (see \cite{Kashi}). For $\mathfrak{gl}(m,n)$ and
$\mathfrak{q}(n),$ it becomes more complicated.\ We postpone the background
useful to prove Proposition \ref{Cor_RSK} to the Appendix.

\section{Pitman transform on the space of paths}

\label{Sec_Pitman}

\subsection{Paths and random walks in $\mathbb{Z}^{n}$}

Recall that $N$ is the cardinality of $\mathcal{A}^{\diamondsuit}$. Denote by
$B^{\diamondsuit}=\{e_{i},i\in\mathcal{A}^{\diamondsuit}\}$ the standard basis
of $\mathbb{Z}^{N}$. We consider paths in $\mathbb{Z}^{N}$ with steps in
$B^{\diamondsuit}$. Observe there is a straightforward bijection between the
set of such paths of length $\ell$ starting from a fixed point $A\in
\mathbb{Z}^{N}$ and the set of words of length $\ell$ on the alphabet $\mathcal{A}
^{\diamondsuit}.\;$More precisely, the word $x_{1}\cdots x_{\ell}
\in(\mathcal{A}^{\diamondsuit})^{\ell}$ of length $\ell$ corresponds to the
path starting at $A$ whose $k$-th step is the translation by $e_{x_{k}}.$
In the sequel we will identify the paths starting from a fixed point $A$ with
the words of $(\mathcal{A}^{\diamondsuit})^{\ell}$.

For any $i\in\mathcal{A}^{\diamondsuit}$, consider $p_{i}\in]0,1[$ and assume
that $\sum_{i\in\mathcal{A}^{\diamondsuit}}p_{i}=1$; this defines a
probability measure on $\mathcal{A}^{\diamondsuit}$ (or equivalently
$B^{\diamondsuit}$). On the space $(\mathcal{A}^{\diamondsuit})^{\mathbb{N}}$
of sequences on the alphabet $\mathcal{A}^{\diamondsuit}$, endowed with the
$\sigma$-algebra $\mathcal{P}\left(  \left(  \mathcal{A}^{\diamondsuit
}\right)  ^{\otimes\mathbb{N}}\right)  $ generated by the cylinder sets, we
consider the infinite product probability measure $\mathbb{P}=p^{\otimes
\mathbb{N}}$; the random variables $X_{\ell},\ell\geq1$, defined on the
probability space $\left(  (\mathcal{A}^{\diamondsuit})^{\mathbb{N}
},\mathcal{P}\left(  \left(  \mathcal{A}^{\diamondsuit}\right)  ^{\otimes
\mathbb{N}}\right)  ,\mathbb{P}\right)  $ by $X_{\ell} :w=(x_{i})_{i\geq
1}\mapsto x_{\ell},$ are independent and identically distributed with law $p$.
Their mean vector
%
is $\displaystyle\mathbf{m}:=\sum_{i\in\mathcal{A}^{\diamondsuit}}p_{i}e_{i}$.

We denote by $\pi_{\infty,\ell}$ the canonical projection from $(\mathcal{A}
^{\diamondsuit})^{\mathbb{N}}$ onto $(\mathcal{A}^{\diamondsuit})^{\ell}$
defined by
\[
\mathrm{for}\text{ }\mathrm{all}\ w=x_{1}x_{2}\ldots\in(\mathcal{A}
^{\diamondsuit})^{\mathbb{N}},\qquad\pi_{\ell}(w)=w^{(l)}:=x_{1}x_{2}\ldots
x_{\ell}.
\]
We set $\mathcal{W}_{\ell}=\mathrm{wt}\circ\pi_{\infty,\ell},\label{def_Wl}$
that is $\mathcal{W}_{\ell}(w)=\mathrm{wt}\left(  w^{(\ell)}\right)  .$ The
random process $\mathcal{W}=\left(  \mathcal{W}_{\ell}\right)  _{\ell
\geq0}$ is a random walk on $\mathbb{Z}^{N}$ since $\mathcal{W}_{\ell
}=X_{1}+\cdots+X_{\ell}$ (here, we consider that the random variables
$X_{\ell}$ take their values in the set ${\mathcal{B}}^{\diamondsuit}$). In
particular $\mathcal{W}$ is a Markov chain with transition matrix
$\Pi_{\mathcal{W}}$ given by
\begin{equation}
\Pi_{\mathcal{W}}(\alpha,\beta)=\left\{
\begin{array}
[c]{l}
p_{i}\text{ if }\beta-\alpha=e_{i}\text{ with }i\in\mathcal{A}^{\diamondsuit
},\\
0\text{ otherwise.}
\end{array}
\right.  \label{Pi(W)}
\end{equation}

\subsection{Pitman transform of paths}

Our aim is now to define a \emph{Pitman transform} $\mathfrak{P}
^{\diamondsuit},\diamondsuit\in\{\emptyset,h,s\}$, on the set of paths we have
considered in the previous section; this will be a generalization of the
classical Pitman transform in the same spirit as in \cite{BBOC1} and it will be defined in fact on $(\mathcal{A}^{\diamondsuit
})^{\mathbb{N}}$.\ 

For any $w\in(\mathcal{A}^{\diamondsuit})^{\mathbb{N}}$, set $\mathfrak{P}
^{\diamondsuit}(w):=\left(  \mathfrak{P}^{\diamondsuit}(w^{(\ell)})\right)
_{\ell}$. We write $sh(T)$ for the shape of the $\diamondsuit$-tableau $T$; it
belongs to $\mathcal{P}^{\diamondsuit}$. The shape of a sequence $(T_{\ell
})_{\ell}$ of tableaux will be the sequence $\left(sh(T_{\ell})\right)_{\ell}$. We now
set
\[
\mathfrak{P}^{\diamondsuit}(w)=sh\left(P^{\diamondsuit}(w)\right)=sh\left(Q^{\diamondsuit
}(w)\right)\text{ for any }w\in(\mathcal{A}^{\diamondsuit})^{\mathbb{N}}.
\]
\begin{example}
Consider $w=1121231212\cdots$. The path in $\mathbb{Z}^{3}$ associated
to $w$ remains in $\mathcal{P}^{\emptyset}$; we obtain
\[
\begin{tabular}
[c]{c|ccccccccccc}
$\ell$ & $1$ & $2$ & $3$ & $4$ & $5$ & $6$ & $7$ & $8$ & $9$ & $10$ & $\cdots$ 
\\\cline{1-12}
\ \\
$sh\left(P^{\emptyset}(w^{(\ell)})\right)$ & $(1)$ & $(2)$ & $(2,1)$ & $(3,1)$ & $(3,2)$
& $(3,2,1)$ & $(4,2,1)$ & $(4,3,1)$ & $(5,3,1)$ & $(5,4,1)$ & $\cdots$
\ \\
\\\cline{1-12}
\ \\
$sh\left(P^{s}(w^{(\ell)})\right)$ & $(1)$ & $(2)$ & $(3)$ & $(3,1)$ & $(4,1)$ & $(5,1)$
& $(5,2)$ & $(5,3)$ & $(5,3,1)$ & $(6,3,1)$ & $\cdots$
\end{tabular}
\]

\end{example}

\bigskip

\noindent\textbf{Remarks: }\noindent$\mathrm{(i)}$ In \cite{LLP}, for a simple Lie algebra $\mathfrak{g}$,  the
generalized Pitman transform  of $w$
was defined from the crystal structure on the set of words as the weight of
the highest weight vertex (source vertex) of the connected component
containing $w.\;$For $\mathfrak{g}=\mathfrak{gl}(n)$ this definition agrees
with that we have just introduced in terms of the insertion algorithm on
tableaux.\ For $\mathfrak{gl}(m,n)$ and $\mathfrak{q}(n),$ there is also a
crystal structure on the set of words but it is more complicated to describe;
in particular there may exist several highest weight vertices for a given
connected component. It becomes thus easier to define the generalized Pitman
transform with the help of insertions algorithm on tableaux. This is what
we do here, in the same spirit of \cite{OC1}.

\noindent$\mathrm{(ii)}$ In view of the previous example, one sees that the
Pitman transform $\mathcal{P}^{s}$ does not fix the paths contained in
$\mathcal{P}^{s}$ (but $\mathcal{P}^{\emptyset}$ does). This phenomenon can be
explained by special behavior of crystal tensor product when $\diamondsuit=h$
or $s$ (see Lemma \ref{Lem_split} and Remark which follows).

\bigskip

We then consider the random variable $\mathcal{H}_{\ell}^{\diamondsuit
}:=\mathcal{W}_{\ell}\circ\mathfrak{P}^{\diamondsuit}$ defined on the
probability space $\left(  (\mathcal{A}^{\diamondsuit})^{\mathbb{N}
},\mathcal{P}\left(  \left(  \mathcal{A}^{\diamondsuit}\right)  ^{\otimes
\mathbb{N}}\right)  ,\mathbb{P}\right)  $, with values in $\mathcal{P}
^{\diamondsuit}.\;$This yields a stochastic process $\mathcal{H}
^{\diamondsuit}\mathcal{=(H}_{\ell}^{\diamondsuit}\mathcal{)}_{\ell\geq0}$.

\begin{proposition}
\label{Prop_laws} For any $\diamondsuit\in\{\emptyset,h,s\},$ any $\ell
\in\mathbb{N}$ and $\lambda\in P_{+}$, one gets
\[
\mathbb{P}[\mathcal{H}_{\ell}^{\diamondsuit}=\lambda]=f_{\lambda
}^{\diamondsuit}\cdot s_{\lambda}^{\diamondsuit}(p).
\]

\end{proposition}

\begin{proof}
By definition of the random variable $\mathcal{H}_{\ell}^{\diamondsuit},$ we
have
\[
\mathbb{P}[\mathcal{H}_{\ell}^{\diamondsuit}=\lambda]=\sum_{T\text{ tableau of
shape }\lambda}\left(  \sum_{w\in B(T)}p_{w}\right)  =\sum_{T\text{ tableau of
shape }\lambda}\mathbb{P}[B(T)].
\]
By (2) of Theorem \ref{Th_RSK}, we have $\mathbb{P}[B(T)]=s_{\lambda
}^{\diamondsuit}(p)$ and in particular it does not depend on $T$ but only on
$\lambda$.\ By (3) of Theorem \ref{Th_RSK}, we then deduce $\mathbb{P}
\left[\mathcal{H}_{\ell}^{\diamondsuit}=\lambda\right]=f_{\lambda}^{\diamondsuit
}\cdot s_{\lambda}^{\diamondsuit}(p)$.
\end{proof}

\bigskip

We can now state the main result of this section.\ For any $\lambda,\mu
\in\mathcal{P}^{\diamondsuit},$ we write $\delta_{\mu\rightsquigarrow\lambda
}^{\diamondsuit}=1$ when $\lambda$ is obtained by adding one box to $\mu$ as in
(\ref{ten_by_V}), and $\delta_{\mu\rightsquigarrow\lambda
}^{\diamondsuit}=0$ otherwise.

\begin{theorem}
\label{Th_main}The stochastic process $\mathcal{H}^{\diamondsuit}$ is a Markov
chain with transition probabilities
\begin{equation}
\Pi_{\mathcal{H}^{\diamondsuit}}(\mu,\lambda)=\frac{s_{\lambda}^{\diamondsuit
}(p)}{s_{\mu}^{\diamondsuit}(p)}\delta_{\mu\rightsquigarrow\lambda
}^{\diamondsuit}\quad\lambda,\mu\in\mathcal{P}^{\diamondsuit}. \label{Pi(H)}
\end{equation}

\end{theorem}

\begin{proof}
Consider a sequence of dominant weights $\lambda^{(1)},\ldots,\lambda^{(\ell
)},\lambda^{(\ell+1)}$ in $\mathcal{P}^{\diamondsuit}$ such that
$\lambda^{(k)}\leadsto\lambda^{(k+1)}$ for any $k=1,\ldots,\ell.\;$We have
seen in \S \ \ref{subsec_RSK} that this determines a unique standard tableau
$T$ and we have
\[
\mathbb{P}[\mathcal{H}_{\ell+1}^{\diamondsuit}=\lambda^{\ell+1},\mathcal{H}
_{k}^{\diamondsuit}=\lambda^{(k)}\text{ for any }k=1,\ldots,\ell]=\sum_{w\in
B^{\diamondsuit}(T)}p_{w}=s_{\lambda^{(\ell+1)}}^{\diamondsuit}(p).
\]
Similarly, we have
\[
\mathbb{P}[\mathcal{H}_{k}^{\diamondsuit}=\lambda^{(k)}\text{ for any
}k=1,\ldots,\ell]=s_{\lambda^{(\ell)}}^{\diamondsuit}(p).
\]
Hence
\[
\mathbb{P}[\mathcal{H}_{\ell+1}^{\diamondsuit}=\lambda^{(\ell+1)}
\mid\mathcal{H}_{k}^{\diamondsuit}=\lambda^{(k)}\text{ for any }
k=1,\ldots,\ell]=\frac{s_{\lambda^{(\ell+1)}}^{\diamondsuit}(p)}
{s_{\lambda^{(\ell)}}^{\diamondsuit}(p)}.
\]
In particular, $\mathbb{P}[\mathcal{H}_{\ell+1}^{\diamondsuit}=\lambda
^{(\ell+1)}\mid\mathcal{H}_{k}^{\diamondsuit}=\lambda^{(k)}$ for any
$k=1,\ldots,\ell]$ depends only on $\lambda^{(\ell+1)}$ and $\mu
=\lambda^{(\ell)}$, this is the Markov property.
\end{proof}

\bigskip

\noindent\textbf{Remarks:} \noindent$\mathrm{(i)}$ By Proposition
\ref{Cor_RSK}, for any $\lambda,\mu$ in $\mathcal{P}^{\diamondsuit}$ with
$\mu\subset\lambda$, the multiplicity $f_{\lambda/\mu}^{\diamondsuit}$ is
equal to the number of standard tableaux of shape $\lambda/\mu$.\ The datum of
such a tableau $T$ is equivalent to that of the sequence $(\mu,\lambda
^{(1)},\ldots,\lambda^{(\ell)})$ of dominant weights where $\ell=\left|
\lambda\right|  -\left|  \mu\right|  $; furthermore, for $k=1, \ldots, \ell$,
the shape of the $\diamondsuit$-diagram $\lambda^{(k)}$ is obtained by adding
to $\mu$ the boxes of $T$ filled by the letters in $\{1,\ldots,k\}.\;$
Therefore there is a bijection between the standard tableaux of shape
$\lambda/\mu$ and the paths from $\mu$ to $\lambda$ which remain in
$\mathcal{P}^{\diamondsuit}$.

\noindent$\mathrm{(ii)}$ Let $\lambda\in\mathcal{P}^{\diamondsuit}$ and
consider the irreducible highest weight representation $V^{\diamondsuit
}(\lambda).\;$Assume $\left\vert \lambda\right\vert =\ell.\;$Let $T$ be a
standard tableau of shape $\lambda$ and define $B^{\diamondsuit}(T)$ as in
Theorem \ref{Th_RSK}.\ Then for any weight $\mu,$ the dimension $K_{\lambda
,\mu}$ of the $\mu$-weight space in $V(\lambda)$ is equal to the number of
words in $B^{\diamondsuit}(T)$ of weight $\mu$. This follows from the
bijection between $B^{\diamondsuit}(T)$ and $T^{\diamondsuit}(\lambda)$
obtained in Theorem~\ref{Th_RSK}.\ Since we have identified paths and words,
the integer $K_{\lambda,\mu}$ is equal to the number of paths from $0$ to
$\mu$ which remains in $B^{\diamondsuit}(T)$.


Recall that $\mathbf{m}=\mathbb{E}(X)$ is the drift of the random walk defined
in (\ref{def_Wl}).\ One gets $\displaystyle\mathbf{m}=\sum_{i=1}^{n}p_{i}
e_{i}$ for $\diamondsuit=\emptyset, s$ and $\displaystyle\mathbf{m}=\sum
_{i=1}^{m}p_{\overline{i}}e_{i}+\sum_{j=1}^{n}p_{j}e_{j+m}$ for $\diamondsuit
=h$. In the sequel, we will assume that the following condition is satisfied:

\begin{condition}
\label{Cond}

\begin{enumerate}
\item For $\diamondsuit=\emptyset,s$ one assumes $p_{1}>\cdots>p_{n}>0$,

\item For $\diamondsuit=h$, one assumes $p_{\overline{m}}>\cdots
>p_{\overline{1}}>0$ and $p_{1}>\cdots>p_{n}>0$.
\end{enumerate}
\end{condition}

%

We will need the following result in Section \ref{Sec_restric}.

\begin{proposition}
\label{Prop_dec_skew}

\begin{enumerate}
\item For any $\lambda,\nu\in\mathcal{P}^{\diamondsuit}$ such that $\nu
\subset\lambda$, one has $\displaystyle f_{\lambda/\nu}^{\diamondsuit}
=\sum_{\mu\in\mathcal{P}^{\diamondsuit}}f_{\mu}^{\diamondsuit}\cdot m_{\mu,\nu
}^{\lambda,\diamondsuit}.$

\item Assume that $\mathbf{m}$ satisfies Condition \ref{Cond}.\ Consider a
sequence of weights $(\lambda^{(a)})_{a\in\mathbb{N}}$ of the form
$\lambda^{(a)}=a\mathbf{m}+o(a)$ and fix a nonnegative integer $\ell$. Then,
for large enough $a$, the weight $\lambda^{(a)}$ belongs to
$\mathcal{P}^{\diamondsuit}$. Moreover, for any $\kappa,\mu\in\mathcal{P}
^{\diamondsuit}$ such that $\left\vert \mu\right\vert =\ell$ and $\left\vert
\lambda^{(a)}\right\vert =\left\vert \kappa\right\vert +\ell$, we have
$m_{\kappa,\mu}^{\lambda^{(a)},\diamondsuit}=K_{\mu,\lambda^{(a)}-\kappa
}^{\diamondsuit}$.

\item For $(\lambda^{(a)})_{a\in\mathbb{N}}$ as above and any
$\mu\in\mathcal{P}^{\diamondsuit}$, one gets, for any large enough $a$,
\begin{equation}
f_{\lambda^{(a)}/\mu}^{\diamondsuit}=\sum_{\kappa\in\mathcal{P}^{\diamondsuit
}}f_{\kappa}^{\diamondsuit}\cdot K_{\mu,\lambda^{(a)}-\kappa}^{\diamondsuit}
=\sum_{\gamma\in P}f_{\lambda^{(a)}-\gamma}^{\diamondsuit}\cdot K_{\mu,\gamma
}^{\diamondsuit}. \label{dec_skew_coef}
\end{equation}

\end{enumerate}
\end{proposition}

\begin{proof}
To prove 1, write $L=\left\vert \lambda\right\vert -\left\vert
\nu\right\vert $. Then, by the definition of the Schur functions, one gets,
using (\ref{skew_multi}) and (\ref{LRcoef}),
\[
s_{\nu}^{\diamondsuit}(s)^{L}=\sum_{\mu}f_{\mu}^{\diamondsuit}\cdot s_{\mu
}^{\diamondsuit}\cdot s_{\nu}^{\diamondsuit}=\sum_{\mu}\sum_{\lambda}f_{\mu
}^{\diamondsuit}\cdot m_{\mu,\nu}^{\lambda,\diamondsuit}\cdot s_{\lambda}^{\diamondsuit
}=\sum_{\lambda}f_{\lambda/\nu}^{\diamondsuit}\cdot s_{\lambda}^{\diamondsuit}
\]
where $s:=s_{(1)}^{\diamondsuit}$ is the character of the natural
representation $V$ and all the sums run over $\mathcal{P}^{\diamondsuit}.$ The
assertion immediately follows by comparing the two last expressions.

To prove 2, observe first that $\lambda^{(a)}$ belongs to $\mathcal{P}
^{\diamondsuit}$ for $a$ sufficiently large because $\mathbf m$ satisfies Condition
\ref{Cond}.\ For any $\kappa\in\mathcal{P}^{\diamondsuit}$ such that
$\left\vert \kappa\right\vert =\left\vert \lambda^{(a)}-\right\vert \ell$, we
have by assertion 1
\[
f_{\lambda^{(a)}/\kappa}^{\diamondsuit}=\sum_{\mu\in\mathcal{P}^{\diamondsuit
},\left\vert \mu\right\vert =\ell}f_{\mu}^{\diamondsuit}\cdot m_{\mu,\kappa
}^{\lambda^{(a)},\diamondsuit}=\sum_{\mu\in\mathcal{P}^{\diamondsuit
},\left\vert \mu\right\vert =\ell}f_{\mu}^{\diamondsuit}\cdot m_{\kappa,\mu
}^{\lambda^{(a)},\diamondsuit}
\]
since $m_{\mu,\kappa}^{\lambda^{(a)},\diamondsuit}=m_{\kappa,\mu}
^{\lambda^{(a)},\diamondsuit}$. Write $K_{\otimes\ell,\lambda^{(a)}-\kappa
}^{\diamondsuit}$ for the dimension of the weight space $\lambda^{(a)}-\kappa$
in $(V^{\diamondsuit})^{\otimes\ell}$. We have
\[
K_{\otimes\ell,\lambda^{(a)}-\kappa}^{\diamondsuit}=\sum_{\mu\in
\mathcal{P}^{\diamondsuit},\left\vert \mu\right\vert =\ell}f_{\mu
}^{\diamondsuit}\cdot K_{\mu,\lambda^{(a)}-\kappa}^{\diamondsuit}
\]
by decomposing $(V^{\diamondsuit})^{\otimes\ell}$ in its irreducible
components.\ By 2 of Proposition \ref{Cor_RSK}, one has $0\leq m_{\kappa
,\mu}^{\lambda^{(a)},\diamondsuit}\leq K_{\mu,\lambda^{(a)}-\kappa
}^{\diamondsuit}$ for any $\mu\in\mathcal{P}^{\diamondsuit}$. It thus suffices
to show that $f_{\lambda^{(a)}/\kappa}^{\diamondsuit}=K_{\otimes\ell
,\lambda^{(a)}-\kappa}^{\diamondsuit}$ for $a$ large enough. Observe that
$K_{\otimes\ell,\lambda^{(a)}-\kappa}^{\diamondsuit}$ is equal to the number
of words of length $\ell$ and weight $\lambda^{(a)}-\kappa$ on $\mathcal{A}
^{\diamondsuit}$. On the other hand, by 1 of Proposition \ref{Cor_RSK}, we
deduce that $f_{\lambda^{(a)}/\kappa}^{\diamondsuit}$ is the number of
sequence of diagrams $(\delta^{(0)},\ldots,\delta^{(\ell)})$ in $\mathcal{P}
^{\diamondsuit}$ of length $\ell$ such that $\delta^{(0)}=\kappa,\delta
^{(\ell)}=\lambda^{(a)}$ and $\delta^{(k+1)}/\delta^{(k)}=(1)$ for any
$k=0,\ldots,\ell-1$.\ For $\mathfrak{g}=\mathfrak{gl}(n)$ or $\mathfrak{q}
(n),$ we associate to $(\delta^{(0)},\ldots,\delta^{(\ell)})$ the word
$w=x_{1}\cdots x_{\ell}$ where for any $k=1,\ldots,\ell-1$, one sets $x_{k}=i$
if the box $\delta^{(k)}/\delta^{(k-1)}$ appears in the row $i\in
\{1,\ldots,n\}$ of $\delta^{(k)}$. For $\mathfrak{g}=\mathfrak{gl}(m,n),$ we
associate similarly to $(\delta^{(0)},\ldots,\delta^{(\ell)})$ the word
$w=x_{1}\cdots x_{\ell}$ where for any $k=1,\ldots,\ell-1$, one sets
$x_{k}=\overline{i}$ if the box $\delta^{(k)}/\delta^{(k-1)}$ appears in the
row $i\in\{1,\ldots,m\}$ of $\delta^{(k)}$ and $x_{k}=j$ if this box appears
in the column $j\in\{1,\ldots,n\}$ of $\delta^{(k)}$. In both cases, one
verifies that this map is injective map from the set of
sequences of diagrams $(\delta^{(0)},\ldots,\delta^{(\ell)})$ into the set of
words of length $\ell$ and weight $\lambda^{(a)}-\kappa$ on $\mathcal{A}
^{\diamondsuit}$. Moreover this map is surjective; indeed we have $\lambda
^{(a)}=a\mathbf{m}+o(a)$ and $\mu$ is fixed, in particular, $\kappa$ and
$\lambda$ differ by at most $\left\vert \mu\right\vert $ boxes. So, by
Condition \ref{Cond}, when $a$ is sufficiently large, adding $\left\vert
\mu\right\vert $ boxes in any order to the rows of $\kappa$ always yields a
diagram in $\mathcal{P}^{\diamondsuit}$. Therefore $f_{\lambda^{(a)}/\kappa
}^{\diamondsuit}=K_{\otimes\ell,\lambda^{(a)}-\kappa}^{\diamondsuit}$.
\end{proof}

\section{Conditioning to stay in $\mathcal{P}^{\diamondsuit}$}

\label{Sec_restric}Denote by $\mathcal{W}^{\diamondsuit}$ the random walk
$\mathcal{W}$ conditioned to stay in $\mathcal{P}^{\diamondsuit}.\;$The aim of
this section is to determine the law of $\mathcal{W}^{\diamondsuit}$. As we
saw in \S \ \ref{mcc}, the process $\mathcal{W}^{\diamondsuit}$ is a Markov
chain; we denote by $\Pi_{\mathcal{W}^{\diamondsuit}}$ its transition matrix.

We have explicit formulae for the transition matrices $\Pi_{\mathcal{H}^{\diamondsuit}}$ and
$\Pi_{\mathcal{W}}$ of the Markov chains $\mathcal{W}$ and
$\mathcal{H}^{\diamondsuit}$. It also follows from (\ref{Pi(H)}) and
(\ref{Pi(W)}) that $\Pi_{\mathcal{H}^{\diamondsuit}}$ is the Doob $\psi
$-transform of $\Pi_{\mathcal{W}}$ where $\psi$ is the harmonic function
\begin{equation}
\psi:\left\{
\begin{array}
[c]{l}
\mathcal{P}^{\diamondsuit}\rightarrow\mathbb{R}_{>0}\\
\lambda\mapsto p^{-\lambda}s_{\lambda}^{\diamondsuit}(p)
\end{array}
\right. . \label{psi}
\end{equation}
On the other hand, we know by \S \ \ref{sub_sec_Doobh} that $\Pi
_{\mathcal{W}^{\diamondsuit}}$ is the Doob $h_{\mathcal{P}^{\diamondsuit}}
$-transform of $\Pi_{\mathcal{W}}$ where $h_{\mathcal{P}^{\diamondsuit}
}(\lambda):=\mathbb{P}\left[  \forall\ell\geq1,\,\mathcal{W}_{\ell}
\in\mathcal{P}^{\diamondsuit}\mid\mathcal{W}_{1}=\lambda\right]  $ for any
$\lambda\in\mathcal{P}^{\diamondsuit}$. We are going to prove that $\psi$ and
$h_{\mathcal{P}^{\diamondsuit}}$ coincide.

\subsection{Limit of $\psi$ along a drift}

Assume Condition \ref{Cond} is satisfied. Then the products
\begin{gather}
\nabla^{\emptyset}=\frac{1}{\prod_{1\leq i<j\leq n}(1-\frac{p_{j}}{p_{i}}
)},\quad\nabla^{h}=\frac{\prod_{i=1}^{m}\prod_{j=1}^{n}(1+\frac{p_{j}
}{p_{\overline{i}}})}{\prod_{\overline{m}\leq\overline{i}<\overline{j}
\leq\overline{1}}(1-\frac{p_{\overline{j}}}{p_{\overline{i}}})\prod_{1\leq
r<s\leq n}(1-\frac{p_{s}}{p_{r}})}\label{def_nabla}\\\text{ and}\quad
\nabla^{s}=\prod_{1\leq i<j\leq n}\dfrac{p_{i}+p_{j}}{p_{i}-p_{j}}\nonumber
\end{gather}
are well-defined.

\begin{proposition}
\label{prop-_imit-psi}Assume Condition \ref{Cond} is satisfied and consider a
sequence $(\lambda^{(a)})_{a\in\mathbb{N}}$ in $\mathcal{P}^{\diamondsuit}$
such that $\lambda^{(a)}=a\mathbf{m}+o(a).$ Then $\lim_{a\rightarrow+\infty
}p^{-\lambda^{(a)}}s_{\lambda^{(a)}}^{\diamondsuit}(p)=\nabla^{\diamondsuit}.$
\end{proposition}

\begin{proof}
Assume first $\diamondsuit=\emptyset$.\ By the Weyl character formula, we
have
\[
s_{\lambda^{(a)}}^{\emptyset}(p)=\nabla^{\emptyset}\sum_{\sigma\in S_{n}
}\varepsilon(\sigma)p^{\sigma(\lambda^{(a)}+\rho)-\rho}.
\]
This gives
\[
p^{-\lambda^{(a)}}s_{\lambda^{(a)}}^{\emptyset}(p)=\nabla^{\emptyset}
\sum_{\sigma\in S_{n}}\varepsilon(\sigma)p^{\sigma(\lambda^{(a)}+\rho
)-\lambda^{(a)}-\rho}.
\]
For $\sigma=1$, one gets $\varepsilon(\sigma)p^{\sigma(\lambda^{(a)}
+\rho)-\lambda^{(a)}-\rho}=1$. So it suffices to prove that
\[
\lim_{a\rightarrow+\infty}\varepsilon(\sigma)p^{\sigma(\lambda^{(a)}
+\rho)-\lambda^{(a)}-\rho}=0
\]
for any $\sigma\neq1$. Consider $\sigma\neq1$ and observe that
\[
\lambda^{(a)}+\rho-\sigma(\lambda^{(a)}+\rho)=\lambda^{(a)}-\sigma
(\lambda^{(a)})+\rho-\sigma(\rho)=a(\mathbf{m}-\sigma(\mathbf{m}))+\rho
-\sigma(\rho)+o(a).
\]
Since $\mathbf{m}$ satisfies Condition \ref{Cond}, the coordinates of
$\mathbf{m}$ strictly decrease and are positive; this implies that
$\displaystyle p^{\mathbf{m}-\sigma(\mathbf{m})}>1$.

Assume $\diamondsuit=h.\;$Since \ref{Cond} is satisfied, for $a$ large enough
the Young diagram of $\lambda^{(a)}$ has a box at position $(i,j)$ for any
$i\in\{1,\ldots,m\}$ and any $j\in\{1,\ldots,n\}$. For such an integer $a,$ we
thus have by (\ref{superWeyl})
\[
s_{\lambda}^{h}(p)=\nabla^{h}\sum_{w\in S_{m}\times S_{n}}\varepsilon
(w)p^{w(\lambda+\rho_{+})-\rho_{+}}.
\]
The arguments are then the same as in the proof of the case $\diamondsuit
=\emptyset$.

Finally, assume $\diamondsuit=s$. Since $\lambda^{(a)}$ has only positive
coordinates, for $a$ large enough, we have by (\ref{rela}) and the case
$\diamondsuit=\emptyset$ (for the sequence $\lambda^{(a)}-\rho$)
\[
\lim_{a\rightarrow+\infty}p^{-\lambda^{(a)}}s_{\lambda^{(a)}}^{s}
(p)=\lim_{a\rightarrow+\infty}p^{-\lambda^{(a)}+\rho}s_{\lambda^{(a)}-\rho
}^{\emptyset}(p)\left(  \prod_{1\leq i<j\leq n}p_{i}+p_{j}\right)  p^{-\rho
}=\nabla^{s}.
\]

\end{proof}

\subsection{The transition matrix $\Pi_{\mathcal{W}^{\diamondsuit}}$}\label{tmP}

In this paragraph, we also assume Condition \ref{Cond} is satisfied.\ Write
$\Pi^{\diamondsuit}$ for the restriction of $\Pi_{\mathcal{W}}$ to
$\mathcal{P}^{\diamondsuit}$. Let us denote by $\Gamma$ the Green function
associated to the substochastic matrix $\Pi^{\diamondsuit}$. For any $\mu, \lambda\in P^{\diamondsuit}$, we
have
\[
\Gamma(\mu,\lambda)=
\begin{cases}
(\Pi^{\diamondsuit})^{\ell}(\mu,\lambda) & \text{ if } \ell=\left\vert
\lambda\right\vert -\left\vert \mu\right\vert \geq0\\
0 & \text{ if } \left\vert \lambda\right\vert <\left\vert \mu\right\vert .
\end{cases}
\]
In particular $\Gamma(\mu,\lambda)=0$ if $\lambda\notin\mathcal{P}
^{\diamondsuit}$. We consider the Martin kernel $K(\mu,\lambda)=\frac
{\Gamma(\mu,\lambda)}{\Gamma(0,\lambda)}.$

In order to apply Theorem \ref{Th_Doob}, we want to prove that, almost surely,
$K(\cdot,\mathcal{W}_{\ell})$ converges everywhere to the harmonic function
$\psi(\cdot)$ defined in (\ref{psi}). By definition of $\Pi^{\diamondsuit}$,
we have
\[
(\Pi^{\diamondsuit})^{\ell}(\mu,\lambda)=\mathrm{card}(ST^{\diamondsuit}(\lambda
/\mu))p^{\lambda-\mu}.
\]
where $ST^{\diamondsuit}(\lambda/\mu)$ is the set of standard $\diamondsuit$-tableaux of
shape $\lambda/\mu$ (see the remark after Theorem~\ref{Th_main}).\ Indeed, all
the paths from $\mu$ to $\lambda$ have the same probability $p^{\lambda-\mu}$
and we have seen there are $\mathrm{card}(ST^{\diamondsuit}(\lambda/\mu))$ such tableaux.\ By
Proposition \ref{Cor_RSK} we have
\[
\mathrm{card}(ST^{\diamondsuit}(\lambda/\mu))=f_{\lambda/\mu}^{\diamondsuit}.
\]
According to Proposition \ref{Prop_dec_skew}, given any sequence
$\lambda^{(a)}$ of weights of the form $\lambda^{(a)}=a\mathbf{m}+o(a)$, we
can write for $a$ large enough
\[
\Gamma(\mu,\lambda^{(a)})=f_{\lambda^{(a)}/\mu}^{\diamondsuit}\cdot p^{\lambda
^{(a)}-\mu}=p^{\lambda^{(a)}-\mu}\sum_{\gamma\in P}f_{\lambda^{(a)}-\gamma
}^{\diamondsuit}\cdot K_{\mu,\gamma}^{\diamondsuit}
\]
Since $\displaystyle
\Gamma(0,\lambda^{(a)})=f_{\lambda^{(a)}}^{\diamondsuit}p^{\lambda^{(a)}},$
this yields for $a$ large enough
\[
K(\mu,\lambda^{(a)})=p^{-\mu}\sum_{\gamma\in P}K_{\mu,\gamma}^{\diamondsuit
}\frac{f_{\lambda^{(a)}-\gamma}^{\diamondsuit}}{f_{\lambda^{(a)}
}^{\diamondsuit}}=p^{-\mu}\sum_{\gamma\in P}K_{\mu,\gamma}^{\diamondsuit
}p^{\gamma}\frac{f_{\lambda^{(a)}-\gamma}^{\diamondsuit}\cdot p^{\lambda
^{(a)}-\gamma}}{f_{\lambda^{(a)}}^{\diamondsuit}\cdot p^{\lambda^{(a)}}}.
\]
Thus
\begin{equation}
K(\mu,\lambda^{(a)})=p^{-\mu}\sum_{\gamma\text{ weight of }V^{\diamondsuit
}(\mu)}K_{\mu,\gamma}^{\diamondsuit}p^{\gamma}\frac{\Gamma(0,\lambda
^{(a)}-\gamma)}{\Gamma(0,\lambda^{(a)})}. \label{expand_K}
\end{equation}
Now we have the following proposition

\begin{proposition}
\label{Prop_quotient}Assume Condition \ref{Cond} is satisfied and consider
sequences $(\lambda^{(a)})_{a\in\mathbb{N}}$ and $(\mu^{(a)})_{a\in\mathbb{N}
}$ in $\mathcal{P}^{\diamondsuit}$ such that $\lambda^{(a)}=a\mathbf{m}
+o(a^{\delta+{\frac{1}{2}}})$ for some $\delta>0$ and $\mu^{(a)}=o(a^{1/2}).$
Then
\[
\lim_{a\rightarrow+\infty}\frac{\Gamma(0,\lambda^{(a)}-\mu^{(a)})}
{\Gamma(0,\lambda^{(a)})}=1.
\]

\end{proposition}

\begin{proof}
We define the semigroup $\mathcal{C}^{\diamondsuit}$ by the one of following
formulae
\begin{align*}
\bullet\quad\mathcal{C}^{\emptyset}  &  =\left\{  (x_{1},x_{2},\ldots
,x_{n})\in\mathbb{R}^{n} \mid x_{1}\geq x_{2}\geq\ldots\geq x_{n}
\geq0\right\}  ,\\
\bullet\quad\mathcal{C}^{s}  &  =\left\{  (x_{1},x_{2},\ldots,x_{n}
)\in\mathcal{C}^{\emptyset}\mid x_{i+1}\neq x_{i}\ \text{if}\ x_{i}
\neq0\right\}  ,\\
\bullet\quad\mathcal{C}^{h}  &  =\big\{(x_{\overline m},x_{\overline{m-1}
},\ldots,x_{\bar1},x_{1},x_{2},\ldots,x_{n})\in\mathbb{R}^{m+n}\\
\  &  \quad\qquad\qquad x_{\overline m}\geq x_{\overline{m-1}}\geq\ldots\geq
x_{\bar1}\geq0,\ x_{1}\geq x_{2}\geq\ldots\geq x_{n}\geq0\ \text{and}
\ \forall i>x_{\bar1}, \,x_{i}=0\big\}.
\end{align*}

In any case, $\mathcal{P}^{\diamondsuit}$ is the intersection of
$\mathcal{C}^{\diamondsuit}$ with the integer lattice and the semigroup
$\mathcal{C}^{\diamondsuit}$ satisfies hypothesis (h1) of \S \ \ref{stayC}.
Moreover the random walk $\mathcal{W}$ satisfies (h3) and, under Condition
\ref{Cond}, it satisfies (h2).

The Proposition \ref{Prop_quotient} is thus a direct consequence of Theorem
\ref{RQT}.
\end{proof}

The strong law of large numbers for square integrable i.i.d. random variables
tells us that for any $\delta>0$, almost surely, one gets
\[
\mathcal{W}_{a}=a\mathbf{m}+o(a^{\delta+\frac12}).
\]
This leads to the following

\begin{corollary}
\label{coro_quotient} Assume Condition \ref{Cond} is satisfied and consider
$\gamma\in P$ a fixed weight. Then
\[
\lim_{a\rightarrow+\infty}\frac{\Gamma(0,\mathcal{W}_{a}-\gamma)}
{\Gamma(0,\mathcal{W}_{a})}=1\text{ (a.s.).}
\]

\end{corollary}

We now may state the main result of this section: we denote by $\left(\mathcal{W}_{\ell
}^{\diamondsuit}\right)_{\ell\geq0}$ the random walk $(\mathcal{W}_{\ell
})_{\ell\geq0}$ conditioned to never exit $\mathcal P^{\diamondsuit}$.

\begin{theorem}
\label{Th_coincide}Assume Condition \ref{Cond} is satisfied. \\The
 Markov chains $(\mathcal{H}_{\ell}^{\diamondsuit
})_{\ell\geq0}$ and $\left(\mathcal{W}_{\ell
}^{\diamondsuit}\right)_{\ell\geq0}$ have the same transition matrix.
\end{theorem}

\begin{proof}
The strong law of large numbers states that $\mathcal{W}_{a}=a\mathbf{m}+o(a)$
almost surely. With (\ref{expand_K}) this implies that almost surely, for $a$
large enough
\[
K(\mu,\mathcal{W}_{a})=p^{-\mu}\sum_{\gamma\text{ weight of }V^{\diamondsuit
}(\mu)}p^{\gamma}\cdot K_{\mu,\gamma}^{\diamondsuit}\frac{\Gamma(0,\mathcal{W}
_{a}-\gamma)}{\Gamma(0,\mathcal{W}_{a})}\text{.}
\]
The weights $\gamma$ run over the set of weights of $V^{\diamondsuit}(\mu)$,
which is finite. Corollary \ref{coro_quotient} then gives
\begin{equation}
L:=\lim_{a\rightarrow+\infty}K(\mu,\mathcal{W}_{a})=p^{-\mu}\sum_{\gamma\text{
weight of }V^{\diamondsuit}(\mu)}p^{\gamma K_{\mu,\gamma}^{\diamondsuit}
}=p^{-\mu}s_{\mu}^{\diamondsuit}(p)=\psi^{\diamondsuit}(\mu)\text{ (a.s.),}
\label{K_tends_psi}
\end{equation}
that is, $L$ coincides with the harmonic function (\ref{psi}). By Theorem
\ref{Th_Doob}, there exists a constant $c$ such that $\psi^{\diamondsuit
}=ch_{\mathcal{P}^{\diamondsuit}}$ where $h_{\mathcal{P}^{\diamondsuit}}$ is
the harmonic function defined in \S \ \ref{sub_sec_Doobh} associated to the
restriction of $(\mathcal{W}_{\ell})_{\ell\geq0}$ to $\mathcal{P}
^{\diamondsuit}$.\ By Theorem \ref{Th_main}, we thus derive
\[
\Pi_{h_{\mathcal{P}^{\diamondsuit}}}(\mu,\lambda)=\Pi_{\psi}(\mu,\lambda
)=\Pi_{\mathcal{H}^{\diamondsuit}}(\mu,\lambda)=\frac{s^{\diamondsuit}_{\lambda}(p)}{s^{\diamondsuit}_{\mu
}(p)}\delta^{\diamondsuit}_{\mu\leadsto\lambda}.
\]

\end{proof}

\subsection{Fewer consequences}

As a direct consequence of Theorem \ref{Th_coincide}, one may state the

\begin{corollary}
\label{yrester} \label{Cor_StayinC} Under the assumptions of Theorem
\ref{Th_coincide}, for any $\lambda\in\mathcal{P}^{\diamondsuit}$, we have
\[
\mathbb{P}_{\lambda}[\mathcal{W}_{\ell}\in\mathcal{P}^{\diamondsuit}\text{ for
all }\ell\geq1]=\frac{p^{-\lambda}s_{\lambda}^{\diamondsuit}(p)}
{\nabla^{\diamondsuit}}
\]
where $\nabla^{\diamondsuit}$ was defined in (\ref{def_nabla}).
\end{corollary}

\begin{proof}
Recall that the function $h_{\mathcal{P}^{\diamondsuit}}:\lambda
\longmapsto\mathbb{P}_{\lambda}\left[\mathcal{W}_{\ell}\in\mathcal{P}
^{\diamondsuit}\text{ for all }\ell\geq1\right]$ is harmonic. By Theorem~
\ref{Th_coincide}, there is a positive constant $c$ such that $\mathbb{P}
_{\lambda}\left[\mathcal{W}_{\ell}\in\mathcal{P}^{\diamondsuit}\text{ for all }\ell
\geq1\right]=cp^{-\lambda}s_{\lambda}^{\diamondsuit}(p)$. Now, for any sequence of
dominant weights $(\lambda^{(a)})_{a}$ such that $\lambda^{(a)}=a\mathbf{m}
+o(a)$, one gets
\[
\lim_{a\rightarrow+\infty}\mathbb{P}_{\lambda^{(a)}}\left[\mathcal{W}_{\ell}
\in\mathcal{P}^{\diamondsuit}\text{ for all }\ell\geq1\right]=\lim_{a\rightarrow
+\infty}\mathbb{P}_{0}\left[\mathcal{W}_{\ell}+\lambda^{(a)}\in\mathcal{P}
^{\diamondsuit}\text{ for all }\ell\geq1\right]=1.
\]
On the other hand, we know by Proposition \ref{prop-_imit-psi} that
$\lim_{a\rightarrow+\infty}p^{-\lambda(a)}s_{\lambda(a)}(p)=\nabla^{\diamondsuit}$.
Therefore $c=\frac{1}{\nabla^{\diamondsuit}}$ and we are done.
\end{proof}

\bigskip

\noindent\textbf{Remark: }Define the open cone $C^\circ:=\{x=(x_{1}
,\ldots,x_{n})\in\mathbb{R}^{n}\mid x_{1}>\cdots>x_{n}>0\};$ it is the common
interior of $\mathcal{C}^{\emptyset}$ and $\mathcal{C}^{s}$.
By the previous corollary, we recover
\[
\mathbb{P}_{\lambda}[\mathcal{W}_{\ell}\in\mathcal{C}^{\emptyset}\text{ for
all }\ell\geq1]=p^{-\lambda}s_{\lambda}^{\emptyset}(p)\prod_{1\leq i<j\leq
n}\left(  1-\frac{p_{j}}{p_{i}}\right)
\]
as established by O'Connell in \cite{OC2}.\ Now for any $\lambda\in
\mathcal{P}_{n}\cap C^\circ$, that is any strict partition with $n$ positive
parts, we have $\mathbb{P}_{\lambda}[\mathcal{W}_{\ell}\in\mathcal{C}^{s}$ for
all $\ell\geq1]=\mathbb{P}_{\lambda}[\mathcal{W}_{\ell}\in C^\circ$ for all
$\ell\geq1]$. So Corollary~\ref{yrester} also gives the probability to stay in
the open cone $C^\circ$. Namely, we have for any $\lambda\in\mathcal{P}
_{n}\cap C^\circ$
\begin{equation}
\mathbb{P}_{\lambda}[\mathcal{W}_{\ell}\in C^\circ\text{ for all }\ell
\geq1]=p^{-\lambda}\cdot s_{\lambda}^{s}(p)\cdot \frac{\displaystyle\prod_{1\leq i<j\leq
n}\left(  p_{i}-p_{j}\right)  }{\displaystyle\prod_{1\leq i<j\leq n}\left(
p_{i}+p_{j}\right)  }.\label{exp1}
\end{equation}
In fact representation theory of $\mathfrak{q}(n)$ is not needed to obtain
this last probability. Indeed, $\lambda\in\mathcal{P}^{s}\cap C^\circ$ implies
that $\lambda-\rho\in\mathcal{P}^{\emptyset}$. This immediately gives
\[
\mathbb{P}_{\lambda}[\mathcal{W}_{\ell}\in C^\circ\text{ for all }\ell
\geq1]=\mathbb{P}_{\lambda-\rho}\left[  \mathcal{W}_{\ell}\in\mathcal{C}
^{\emptyset}\text{ for all }\ell\geq1\right]
\]
so that
\begin{equation}
\mathbb{P}_{\lambda}[\mathcal{W}_{\ell}\in C^\circ\text{ for all }\ell
\geq1]=p^{\rho}\cdot p^{-\lambda}\cdot s_{\lambda-\rho}^{\emptyset}(p)\prod_{1\leq i<j\leq
n}\left(  1-\frac{p_{j}}{p_{i}}\right)  .\label{exp2}
\end{equation}
and, by (\ref{rela}), the two expressions (\ref{exp1}) and (\ref{exp2})
actually coincide. Nevertheless, (\ref{exp2}) does no longer hold when
$\lambda\in\mathcal{P}^{s}$ but $\lambda\notin C^\circ$ (that is $\lambda$ has
zero coordinates). In this case the probability to stay in $C^\circ$ is only
given by (\ref{exp1}).

\bigskip

We can also recover a result by Stanley exposed in \cite{St} giving the
asymptotic behavior of $f_{\lambda/\mu}^{\emptyset}$ and extend it to obtain
the asymptotic behavior of $f_{\lambda/\mu}^{s}$ and $f_{\lambda/\mu}^{h}$

\begin{theorem}
\label{Th_Asympt}Suppose that the vector $\mathbf{m}$ satisfies Condition
\ref{Cond} and consider $\mu\in\mathcal{P}^{\diamondsuit}$. If $\lambda
^{(\ell)}=\ell\mathbf{m}+o(\ell^{\alpha})$ with $\alpha<2/3$, then
\begin{equation}
\lim_{\ell\rightarrow\infty}\frac{f_{\lambda^{(\ell)}/\mu}^{\diamondsuit}
}{f_{\lambda^{(\ell)}}^{\diamondsuit}}=s_{\mu}^{\diamondsuit}(p).
\label{asympt}
\end{equation}

\end{theorem}

\begin{proof}
Consider $\lambda^{(\ell)}=\ell\mathbf{m}+o(\ell^{\alpha})$ a sequence of
dominant weights. By Proposition \ref{Prop_dec_skew}
\begin{equation}
\frac{f_{\lambda^{(\ell)}/\mu}^{\diamondsuit}}{f_{\lambda^{(\ell)}
}^{\diamondsuit}}=\sum_{\gamma\in P}K_{\mu,\gamma}^{\diamondsuit}
\frac{f_{\lambda^{(\ell)}-\gamma}^{\diamondsuit}}{f_{\lambda^{(\ell)}
}^{\diamondsuit}}=\sum_{\gamma\in P}K_{\mu,\gamma}^{\diamondsuit}
\frac{f_{\lambda^{(\ell)}-\gamma}^{\diamondsuit}\cdot p^{\lambda^{(\ell)}-\gamma}
}{f_{\lambda^{(\ell)}}^{\diamondsuit}\cdot p^{\lambda^{(\ell)}}}p^{\gamma}
\label{quo-f}
\end{equation}
where the sums are finite since the set of weights in $V^{\diamondsuit}(\mu)$
is finite. Note that, for any $\gamma\in P$
\[
\frac{f_{\lambda^{(\ell)}-\gamma}^{\diamondsuit}\cdot p^{\lambda^{(\ell)}-\gamma}
}{f_{\lambda^{(\ell)}}^{\diamondsuit}p^{\lambda^{(\ell)}}}=\frac
{\mathbb{P}[\mathcal{W}_{1}\in\mathcal{P}^{\diamondsuit},\ldots,\mathcal{W}
_{\ell}\in\mathcal{P}^{\diamondsuit},\mathcal{W}_{\ell}=\lambda^{(\ell
)}-\gamma]}{\mathbb{P}[\mathcal{W}_{1}\in\mathcal{P}^{\diamondsuit}
,\ldots,\mathcal{W}_{\ell}\in\mathcal{P}^{\diamondsuit},\mathcal{W}_{\ell
}=\lambda^{(\ell)}]}.
\]
By Theorem \ref{RQT}, we know that this quotient tends to $1$ when $\ell$
tends to infinity. This implies
\[
\lim_{\ell\rightarrow+\infty}\frac{f_{\lambda^{(\ell)}/\mu}^{\diamondsuit}
}{f_{\lambda^{(\ell)}}^{\diamondsuit}}=\sum_{\gamma\in P}K_{\mu,\gamma
}^{\diamondsuit}p^{\gamma}=s_{\mu}^{\diamondsuit}(p).
\]

\end{proof}

\subsection{Conditioning in dimension 2}

Let us illustrate our results by considering the situation in dimension
$2$.\ The relevant algebras are then $\mathfrak{g}=\mathfrak{gl}
(2),\mathfrak{gl}(1,1)$ or $\mathfrak{q}(2)$. In this very simple case, we have

\begin{enumerate}
\item $\mathcal{P}^{\emptyset}=\{\lambda=(\lambda_{1},\lambda_{2}
)\in\mathbb{Z}^{2}\mid\lambda_{1}\geq\lambda_{2}\geq0\},$

\item $\mathcal{P}^{s}=\{\lambda=(\lambda_{1},\lambda_{2})\in\mathcal{P}
^{\emptyset}\mid\lambda_{1}=\lambda_{2}\Longrightarrow$ $\lambda=(0,0)\},$

\item $\mathcal{P}^{h}=\{\lambda=(\lambda_{\overline{1}},\lambda_{1}
)\in\mathbb{Z}_{\geq0}^{2}\mid\lambda_{\overline{1}}=0\Longrightarrow
\lambda_{1}=0\}.$
\end{enumerate}

By Applying Corollary \ref{Cor_StayinC}, we get
\[
\mathbb{P}_{(0,0)}[\mathcal{W}_{\ell}\in\mathcal{P}^{\diamondsuit}\text{ for
all }\ell\geq1]=\left\{
\begin{array}
[c]{l}
1-\frac{p_{2}}{p_{1}}\text{ for }\diamondsuit=\emptyset,\\
p_{1}(1-\frac{p_{2}}{p_{1}})\text{ for }\diamondsuit=s,\\
p_{\overline{1}}\text{ for }\diamondsuit=h.
\end{array}
\right.
\]
These formulas are not surprising and can be obtained in an elementary way by
interpreting the one-way simple walk as a the random walk on $\mathbb{Z}$ with
transitions $\pm1$. Conditioning to stay in $\mathcal{P}^{s}$ starting from
$(0,0)$ is equivalent to impose first that $\mathcal{W}_{1}=e_{1}$, next that
$\mathcal{W}_{\ell},\ell\geq2,$ remains in the translated of $\mathcal{P}
^{\emptyset}$ by $(1,0)$. Conditioning to stay in $\mathcal{P}^{h}$ simply
means that $\mathcal{W}_{1}=e_{\overline{1}}$.\ 

\noindent Observe the reduction of the random walk $(\mathcal{W}_{\ell}
)_{\ell\geq0}$ in dimension $2$ to a random walk on $\mathbb{Z}$ cannot be
generalized to dimensions greater than $2$.

\noindent In dimension $2$, there is also a simple relation between the Pitman
transforms $\mathfrak{P}^{\emptyset}$ and $\mathfrak{P}^{s}$.\ Consider
$w=x_{1}\cdots x_{\ell}x_{\ell+1}\in(\mathcal{A}
^{\diamondsuit})^{\ell+1}$ and set $\mathfrak{P}^{\emptyset}(w^{\flat
})=(\lambda_{1},\lambda_{2})\in\mathcal{P}^{\emptyset}$ where $w^{\flat}
=x_{1}\cdots x_{\ell}$. Then
\begin{equation}
\mathfrak{P}^{s}(w)=(\lambda_{1}+1,\lambda_{2}). \label{Pit2}
\end{equation}
This means that $(\mathcal{H}_{\ell}^{s})_{\ell\geq0}$ is obtained from
$(\mathcal{H}_{\ell}^{\emptyset})_{\ell\geq0}$ by translation by $(1,0)$ in
space and translation by $1$ in time. Let $w^{s}$ and $w^{\emptyset}$ be the
highest weight vertices associated to $w$ for the $\mathfrak{q}(2)$ and
$\mathfrak{gl}(2)$-structures, respectively. Write $t$ for the rightmost
letter in $w^{\emptyset}$. By using crystal basis theory and relations (1.7)
and (1.8) in \cite{GHSKM}, one can prove that $w^{s}=w^{\emptyset}$ if $t=1$.
If $t=2$, the weight $w^{s}$ is obtained from $w^{\emptyset}$ by changing $t=2$ in $1$.
This gives (\ref{Pit2}).

\noindent For $n>2$, there is no simple relation between the generalized
Pitman transforms $\mathfrak{P}^{\emptyset}$ and $\mathfrak{P}^{s}$ as shown
by the following example.

\begin{examples}
In view of example \ref{ex_s}, we have for $w=232145331$
\[
\begin{tabular}
[c]{c|ccccccccc}
$\ell$ & $1$ & $2$ & $3$ & $4$ & $5$ & $6$ & $7$ & $8$ & $9$\\\cline{2-10}
$\mathcal{H}_{\ell}^{\emptyset}$ & $(1)$ & $(1,1)$ & $(2,1)$ & $(3,1)$ &
$(3,1,1)$ & $(3,1,1,1)$ & $(3,2,1,1)$ & $(3,3,1,1)$ & $(4,3,1,1)$\\
$\mathcal{H}_{\ell}^{s}$ & $(1)$ & $(2)$ & $(2,1)$ & $(3,1)$ & $(4,1)$ &
$(5,1)$ & $(5,2)$ & $(5,2,1)$ & $(5,3,1)$
\end{tabular}
\]
where the middle row is obtained by computing $P^{\emptyset}(w)$.
\end{examples}

\section{Appendix}

\label{Sec-append}The aim of this section is to prove Proposition
\ref{Cor_RSK}.\ To do this, we need first to interpret the RSK correspondence
in terms of crystal basis theory.

\subsection{RSK correspondence and crystal basis theory}
We refer the reader to \cite{Kashi} ($\diamondsuit=\emptyset$), \cite{BKK}
($\diamondsuit=h)$ and \cite{GHSKM} ($\diamondsuit=s$) for  detailed
expositions. To each irreducible representation $V^{\diamondsuit}(\lambda)$ is
associated its crystal graph $B^{\diamondsuit}(\lambda)$ which is an oriented
graph with arrows $\overset{i}{\rightarrow}$ colored by symbols $i\in
\{1,\ldots,n-1\}$ for $\mathfrak{g}=\mathfrak{gl}(n),$ $i\in\{\overline
{m-1},\ldots,\overline{1},0,1,\ldots,n-1\}$ for $\mathfrak{g}=\mathfrak{gl}
(m,n)$ and $i\in\{1,\ldots,n-1\}$ for $\mathfrak{g}=\mathfrak{q}(n)$. More
generally, any representation $M^{\diamondsuit}$ (irreducible or not) appearing
in a tensor product $(V^{\diamondsuit})^{\otimes\ell}$ admits a crystal
$B^{\diamondsuit}(M)$. The crystal $B^{\diamondsuit}(M)$ is graded by the
weights of $\mathfrak{g}$. There is a map $\mathrm{wt:}B^{\diamondsuit
}(M)\rightarrow P$.\ To obtain the decomposition of $M^{\diamondsuit}$ in its
irreducible components, it then suffices to obtain the decomposition of
$B^{\diamondsuit}(M)$ into its connected components. The crystal
$B^{\diamondsuit}(M\otimes N)=B^{\diamondsuit}(M)\otimes B^{\diamondsuit}(N)$
associated to the tensor product $M\otimes N$ of two representations can be
constructed from the crystal of $M$ and $N$ by simple combinatorial rules. A
highest weight vertex in $B^{\diamondsuit}(M)$ is a vertex $b$ for which no
arrow $b^{\prime}\overset{i}{\rightarrow}b$ exists in $B^{\diamondsuit}
(M)$; a lowest weight vertex in $B^{\diamondsuit}(M)$ is a vertex $b$ for
which no arrow $b\overset{i}{\rightarrow}b^{\prime}$ exists in
$B^{\diamondsuit}(M)$.

Let $\sigma_{0}$ be the element of the Weyl group of $\mathfrak{g}$ defined by
\[
\left\{
\begin{array}
[c]{l}
\sigma_{0}(\beta_{1}, \beta_{2}, \ldots,\beta_{n})=(\beta_{n}, \beta_{n-1}, \ldots,\beta_{1})\text{ for
any }\beta\in\mathbb{Z}^{n}\text{ when }\diamondsuit=\emptyset,s,\\
\sigma_{0}(\beta_{\overline{m}},\ldots,\beta_{\overline{1}},\beta_{1}
,\ldots,\beta_{n})=(\beta_{\overline{1}},\ldots\beta_{\overline{m}},\beta
_{n},\ldots,\beta_{1})\text{ for any }\beta\in\mathbb{Z}^{m+n}\text{ when
}\diamondsuit=h.
\end{array}
\right.
\]
For $\diamondsuit=\emptyset\mathfrak{\ }$and $\lambda\in\mathcal{P}
^{\emptyset}$, the crystal $B^{\emptyset}(\lambda)$ contains a unique highest
weight vertex $b^{\emptyset,\lambda}$ and a unique lowest weight vertex
$b_{\lambda}^{\emptyset}$. They have respective weights $\lambda$ and
$\sigma_{0}(\lambda)$. For $\diamondsuit=h,s\mathfrak{\ }$and $\lambda
\in\mathcal{P}^{\diamondsuit}$, the crystal $B^{\diamondsuit}(\lambda)$
contains a unique highest weight vertex $b^{\diamondsuit,\lambda}$ and a
unique lowest weight vertex $b_{\lambda}^{\diamondsuit}$ with weights
$\lambda$ and $\sigma_{0}(\lambda)$. Nevertheless, the crystal
$B^{\diamondsuit}(\lambda)$ may admit other highest weight vertices than
$b^{\diamondsuit,\lambda}$ and similarly other lowest weight vertices than
$b_{\lambda}^{\diamondsuit}$, and the weights  of these vertices are thus distinct from $\lambda$
and $\sigma_{0}(\lambda)$.

\noindent By crystal basis theory, for $\diamondsuit=\emptyset,h$ or $s$, the
multiplicity of $V^{\diamondsuit}(\lambda)$ in $M^{\diamondsuit}$ is given by
the number of highest weight vertices of weight $\lambda$ in $B^{\diamondsuit
}(M)$ or equivalently by the number of its lowest weight vertices of weight
$\sigma_{0}(\lambda)$.

\begin{lemma}
\label{Lem_split}Assume $\diamondsuit=\emptyset$ or $\diamondsuit=s$ and
consider $u\otimes v\in B^{\diamondsuit}(M)\otimes B^{\diamondsuit}(N)$. Then

\begin{enumerate}
\item for $\diamondsuit=\emptyset$, the vertex $u\otimes v$ is a highest
weight vertex only if $u$ is a highest weight vertex,

\item for $\diamondsuit=s$, the vertex $u\otimes v$ is a lowest weight vertex
only if $v$ is a lowest weight vertex.
\end{enumerate}
\end{lemma}

\noindent\textbf{Remark: }Assertion 1 does not hold in general for
$\diamondsuit=h$ or $s$. Moreover Assertion 2 also fails when $\diamondsuit
=h$ which causes some complications. The lack of Assertion 1 for
$\diamondsuit=h,s$ explains also why in these cases the paths which remain
in $\mathcal{P}^{\diamondsuit}$ are not fixed by the Pitman transforms we have defined.

\bigskip

We give below the crystal $B^{\diamondsuit}$ of the defining representation
$V^{\diamondsuit}$.
\begin{gather*}
1\overset{1}{\rightarrow}2\overset{2}{\rightarrow}\cdots\overset
{n-1}{\rightarrow}n\text{ for }\mathfrak{g}=\mathfrak{gl}(n),\\
\overline{m}\overset{\overline{m-1}}{\rightarrow}\overline{m-1}\overset
{\overline{m-2}}{\rightarrow}\cdots\overset{\overline{1}}{\rightarrow
}\overline{1}\overset{0}{\rightarrow}1\overset{1}{\rightarrow}2\overset
{2}{\rightarrow}\cdots\overset{n-1}{\rightarrow}n\text{ for }\mathfrak{g}
=\mathfrak{gl}(m,n),\\
1
\genfrac{.}{.}{0pt}{}{\overset{1}{\rightarrow}}{\underset{\overline{1}
}{\dashrightarrow}}
2\overset{2}{\rightarrow}\cdots\overset{n-1}{\rightarrow}n\text{ for
}\mathfrak{g}=\mathfrak{q}(n).
\end{gather*}
Observe that the vertices of $B^{\diamondsuit}$ coincide with the letters of
$\mathcal{A}^{\diamondsuit}$. The vertices of $(B^{\diamondsuit})^{\otimes
\ell}$ are thus labelled by the words of $(\mathcal{A}^{\diamondsuit})^{\ell}$
by identifying each vertex $b=x_{1}\otimes\cdots\otimes x_{\ell}$ of
$(B^{\diamondsuit})^{\otimes\ell}$ with the word $b=x_{1}\cdots x_{\ell}$.

Two crystals $B$ and $B^{\prime}$ are isomorphic when there exists a bijection
$\phi:B\rightarrow B^{\prime}$ which respects the graph structure i.e. such
that $\phi(a)\overset{i}{\rightarrow}\phi(b)$ in $B^{\prime}$ if and only if
$a\overset{i}{\rightarrow}b$ in $B$. When $B$ and $B^{\prime}$ are crystals
associated to irreducible representations, such a crystal isomorphism exists
if and only if these representations are isomorphic and in that case it is
unique. For any $w\in(B^{\diamondsuit})^{\otimes\ell}$, write $B^{\diamondsuit
}(w)$ for the connected component of $(B^{\diamondsuit})^{\otimes\ell}$
containing $w$. We can now interpret Theorem \ref{Th_RSK} in terms of crystal
basis theory.

\begin{theorem}
\label{Th_RSKcrystal}Consider $w_{1}$ and $w_{2}$ two vertices of
$(B^{\diamondsuit})^{\otimes\ell}$. Then

\begin{enumerate}
\item $P^{\diamondsuit}(w_{1})=P^{\diamondsuit}(w_{2})$ if and only if
$B^{\diamondsuit}(w_{1})$ is isomorphic to $B^{\diamondsuit}(w_{2})$ and the
unique associated isomorphism sends $w_{1}$ on $w_{2}.$

\item $Q^{\diamondsuit}(w_{1})=Q^{\diamondsuit}(w_{2})$ if and only if
$B^{\diamondsuit}(w_{1})=B^{\diamondsuit}(w_{2})$.

\item For any standard $\diamondsuit$-tableau $T,$ the set of words
$B^{\diamondsuit}(T)$ defined in \S \ \ref{subsec_RSK} has the structure of a
crystal graph isomorphic to the abstract crystal $B^{\diamondsuit}(\lambda)$
where $\lambda$ is the shape of $T$.

\item If we denote by $\phi:B^{\diamondsuit}(\lambda)\rightarrow
B^{\diamondsuit}(T)$ this isomorphism, we have $\mathrm{wt}(b)=\mathrm{wt}
(\phi(b))$ for any $b\in B^{\diamondsuit}(T),$ that is the weight graduation
defined on the abstract crystal $B^{\diamondsuit}(T)$ is compatible with the
weight graduation defined on words.
\end{enumerate}
\end{theorem}

\subsection{Proof of Proposition \ref{Cor_RSK}}

Consider $\mu\in\mathcal{P}^{\diamondsuit}$ and $T$ a standard tableau of
shape $\mu$. Let $\ell$ be a nonnegative integer and $\mathcal{U}_{\ell
,T}^{\diamondsuit}$ be the set of pairs $(P,Q)$ where $Q$ is a $\diamondsuit
$-standard tableau with $\ell+\left\vert \mu\right\vert $ boxes containing $T$
as a subtableau (that is $T$ is the subtableau obtained by considering only
the entries $1,\ldots,\left\vert \mu\right\vert $ of $Q$) and $P$ a
$\diamondsuit$-tableau with the same shape as $Q$. By Theorem \ref{Th_RSK},
the restriction of $\theta_{\ell+\left\vert \mu\right\vert }^{\diamondsuit}$
to the subset $B^{\diamondsuit}(T)\times(\mathcal{A}^{\diamondsuit})^{\ell
}\subset\mathcal{A}_{\ell+\left\vert \mu\right\vert }^{\diamondsuit}$ yields
the one to one correspondence
\[
\theta_{\ell,T}^{\diamondsuit}:\left\{
\begin{array}
[c]{l}
B^{\diamondsuit}(T)\otimes(B^{\diamondsuit})^{\otimes\ell}\rightarrow
\mathcal{U}_{\ell,T}^{\diamondsuit}\\
w_{T}\otimes w\mapsto(P^{\diamondsuit}(w_{T}w),Q^{\diamondsuit}(w_{T}w))
\end{array}
\right.  .
\]
Indeed, for any $u\in(B^{\diamondsuit})^{\otimes\ell+\left\vert \mu\right\vert
}$, $T$ is a subtableau of $Q^{\diamondsuit}(u)$ if and only if $u$ can be
written $u=w_{T}\otimes w$ with $w_{T}\in B^{\diamondsuit}(T)$ and
$w\in(B^{\diamondsuit})^{\otimes\ell}$. By crystal basis theory, for
$(B^{\diamondsuit})^{\otimes\ell+\left\vert \mu\right\vert }$, the number of
connected components of $B^{\diamondsuit}(T)\otimes(B^{\diamondsuit}
)^{\otimes\ell}$ isomorphic to some $B^{\diamondsuit}(\lambda)$ is equal to
$f_{\lambda/\mu}^{\diamondsuit}$; using Theorem \ref{Th_RSKcrystal}, we see
that it coincides with the number of standard $\diamondsuit$-tableaux of shape
$\lambda$ containing $T$ as a subtableau. There is a natural bijection
between the set of such tableaux and the set of skew standard $\diamondsuit
$-tableaux of shape $\lambda/\mu$; to any standard tableau $Q$ containing $T$,
it associates the skew tableau obtained by deleting the boxes of $T$ in $Q$
and substracting $\left\vert \mu\right\vert $ to the entries of the remaining
boxes. This proves assertion 1 of Proposition \ref{Cor_RSK}.

\bigskip

Our method to prove Assertion 2 of Proposition \ref{Cor_RSK} depends on
$\diamondsuit$. 

For $\diamondsuit=\emptyset,$ let $b=b_{1}\otimes b_{2}$ be a
highest weight vertex in $B^{\emptyset}(\kappa)\otimes B^{\emptyset}(\mu)$ of
weight $\lambda$; by Lemma \ref{Lem_split}, we must have $b_{1}=b^{\emptyset
,\kappa}$.\ Since $\mathrm{wt}(b)=\mathrm{wt}(b_{1})+\mathrm{wt}(b_{2})$, this
implies that $b_{2}\in B^{\diamondsuit}(\mu)$ has weight $\lambda-\kappa
$.\ Therefore $m_{\kappa,\mu}^{\emptyset,\lambda}\leq K_{\mu,\lambda-\kappa
}^{\emptyset}$ since we can define an injective map from the set of highest
weight vertices of weight $\lambda$ in $B^{\emptyset}(\kappa)\otimes
B^{\emptyset}(\mu)$ to the set of vertices of weight $\lambda-\kappa$ in
$B^{\emptyset}(\mu)$.

For $\diamondsuit=s,$ let $b=b_{1}\otimes b_{2}$ be a lowest weight
vertex in $B^{s}(\kappa)\otimes B^{s}(\mu)$ of weight $\lambda$; we must have
$b_{2}=b_{\mu}^{s}$, so $\mathrm{w}(b_{1})=\sigma_{0}(\lambda)-\sigma_{0}
(\mu)=\sigma_{0}(\lambda-\mu)$. One deduces similarly that\ $m_{\kappa,\mu
}^{s,\lambda}=m_{\mu,\kappa}^{s,\lambda}\leq K_{\kappa,\sigma_{0}(\lambda
-\mu)}^{s}$. But $K_{\kappa,\sigma_{0}(\lambda-\mu)}^{s}=K_{\kappa,\lambda
-\mu}^{s}$ since $\sigma_{0}$ belongs to the Weyl group of $\mathfrak{q}(n)$.
Thus $m_{\mu,\kappa}^{s,\lambda}\leq K_{\kappa,\lambda-\mu}^{s}$ as desired.

It remains to consider the case when $\diamondsuit=h$, i.e. $\mathfrak{gl}
(m,n)$. Since Lemma \ref{Lem_split} is no longer true in this case, we shall
need a different strategy and use the Littlewood-Richardson rule established
in \cite{KK}. One says that a word $w=x_{1}\cdots x_{\ell}$ with letters in
$\mathbb{Z}_{>0}$ is a \textit{permutation word} when for any $k=1,\ldots
,\ell$ and any positive integer $i$, the number of letters $i$ in the prefix
$w^{(k)}=x_{1}\cdots x_{k}$ is greater or equal to the number of letters $i+1$.

Consider $\lambda,\mu,\kappa$ in $\mathcal{P}^{h}$ such that $\left\vert
\mu\right\vert =\left\vert \lambda\right\vert -\left\vert \kappa\right\vert $.
Decompose the diagram of $\lambda$ in $\lambda^{(1)}$ (obtained by considering
its first $m$-rows) and $\lambda^{(2)}$ as in
\S \ \ref{subsec_YD}.\ We denote by $LR_{\mu,\kappa}^{\lambda}$ the set of
tableaux $T$ obtained by filling the Young diagram $\lambda/\kappa$ with
positive integers (see Example \ref{exam_teta}) such that

\begin{enumerate}
\item the rows of $T$ weakly increase from left to right,

\item the columns of $T$ strictly increase from top to bottom,

\item the word $w$ obtained by reading first the $m$ rows of $\left(\lambda/\kappa\right)^{(1)}$ from
right to left and top to bottom, next the $n$ columns of $\left(\lambda/\kappa\right)^{(2)}$ from
right to left and top to bottom is a permutation word such that for any
integer $k\geq1$, the number of letters $k$ in $w$ is equal to the length of
the $k$-row of the Young diagram of $\mu$, that is the $k$-th coordinate of
$(\mu^{(1)},(\mu^{(2)})^{\prime})$ where $(\mu^{(2)})^{\prime}$ is
the conjugate partition of $\mu^{(2)}$.
\end{enumerate}

\begin{proposition}
\cite{KK}With the previous notations, we have $m_{\mu,\kappa}^{h,\lambda
}=\mathrm{card}(LR_{\mu,\kappa}^{\lambda})$.
\end{proposition}

\noindent\textbf{Remark: }The result obtained in \cite{KK} is more general. We
can replace the reading of the LR-tableaux used in our definition by any
admissible fixed reading, that is any reading such that a box $b$ is read
before a box $b^{\prime}$ whenever $b$ is located at the north-west of
$b^{\prime}$.\ 

\bigskip

By the previous proposition, in order to prove that $m_{\kappa,\mu}
^{h,\lambda}\leq K_{\mu,\lambda-\kappa}^{h}$ it suffices to construct an
embedding $\theta$ from $LR_{\mu,\kappa}^{\lambda}$ to the set of $h$-tableaux
of shape $\mu$ and weight $\lambda-\kappa$. We proceed as follows. Consider
$T\in LR_{\mu,\kappa}^{\lambda}$ and write $R_{\overline{m}},\ldots
,R_{\overline{1}}$ the rows with boxes in $\lambda^{(1)}$ and $C_{1}
,\ldots,C_{n}$ the columns with boxes in $\lambda^{(2)}$. In particular, for
any $i=1,\ldots,m,$ and $j=1,\ldots,n$, the rows $R_{\overline{i}}$ contains
$\lambda_{\overline{i}}^{(1)}-\kappa_{\overline{i}}^{(1)}$ letters and the
column $C_{j}$ contains $\lambda_{j}^{(2)}-\kappa_{j}^{(2)}$ letters.

We want to define the $h$-tableau $\theta(T)$. We first construct recursively
a tableau $T^{(1)}$ with $\lambda_{k}^{(1)}-\kappa_{k}^{(1)}$ letters
$\overline{k}$ for any $k=1,\ldots,m$.\ First define $T_{\overline{m}}$ as the
row with $\lambda_{\overline{m}}^{(1)}-\kappa_{\overline{m}}^{(1)}$ letters
$\overline{m}$. Assume $T_{\overline{k+1}},k\in\{2,\ldots,m\}$ is constructed.
Then $T_{\overline{k}}$ is obtained by adding a letter $\overline{k}$ on the
$i$-th row of $T_{\overline{k+1}}$ for each occurrence of the integer $i$ in
$R_{\overline{k}}$.\ 

Set $T^{(1)}=T_{\overline{1}}$ and let us construct $\theta(T)$ by adding
successively letters to $T^{(1)}$. First define $T_{1}$ by adding a letter $1$
on the $i$-th row of $T_{\overline{1}}$ for each integer $i$ in $C_{1}$. Now,
when $T_{k}$ is given, for some $k\in\{1,\ldots,m-1\}$, one constructs
$T_{k+1}$ adding a letter $k$ on the $i$-th row of $T_{k}$ for each integer
$i$ appearing in $C_{k}$.\ Finally, set $\theta(T)=T_{n}$. Observe that by
construction, the $k$-th row of $\theta(T)$ contains as many boxes as the
number of letters $k$ in the permutation word $w$ associated to $T$ by the
previous Assertion $3$.\ Moreover, $\theta(T)$ is a $h$-tableau since $w$ is a
permutation word; it has shape $\mu$ and weight $\lambda-\kappa$ as desired.
Moreover the map $\theta$ is injective since $\theta(T)$ records both the
number and the positions of the letters $k$ in skew shape $\lambda-\kappa$.
This proves that $m_{\kappa,\mu}^{h,\lambda}\leq K_{\mu,\lambda-\kappa}^{h}$
as desired.

\begin{example}
\label{exam_teta}Take $\lambda=(3,3,3\mid3,3)$ with $m=n=3,$ $\kappa
=(2,0,0\mid0,0)$ and $\mu=(3,3,2\mid3,2)=(3,3,2,2,2,1)$. For
\[
T=
\begin{tabular}
[c]{ll|l}\cline{3-3}
&  & \multicolumn{1}{|l|}{$1$}\\\hline
\multicolumn{1}{|l}{$1$} & \multicolumn{1}{|l|}{$1$} &
\multicolumn{1}{|l|}{$2$}\\\hline
\multicolumn{1}{|l}{$2$} & \multicolumn{1}{|l|}{$2$} &
\multicolumn{1}{|l|}{$3$}\\\hline
\multicolumn{1}{|l}{$\mathbf{3}$} & \multicolumn{1}{|l|}{$\mathbf{4}$} &
\\\cline{1-2}
\multicolumn{1}{|l}{$\mathbf{4}$} & \multicolumn{1}{|l|}{$\mathbf{5}$} &
\\\cline{1-2}
\multicolumn{1}{|l}{$\mathbf{5}$} & \multicolumn{1}{|l|}{$\mathbf{6}$} &
\\\cline{1-2}
\end{tabular}
\ \ \text{, we get }\theta(T)=
\begin{tabular}
[c]{|l|ll}\hline
$\bar{3}$ & $\bar{2}$ & \multicolumn{1}{|l|}{$\bar{2}$}\\\hline
$\bar{2}$ & $\bar{1}$ & \multicolumn{1}{|l|}{$\bar{1}$}\\\hline
$\bar{1}$ & $\mathbf{1}$ & \multicolumn{1}{|l}{}\\\cline{1-2}\cline{2-2}
$\mathbf{1}$ & $\mathbf{2}$ & \multicolumn{1}{|l}{}\\\cline{1-2}\cline{2-2}
$\mathbf{1}$ & $\mathbf{2}$ & \multicolumn{1}{|l}{}\\\cline{1-2}
$\mathbf{2}$ &  & \\\cline{1-1}
\end{tabular}
\ \
\]
where $w=1211322456345$ is a permutation word.
\end{example}

\end{document}